\documentclass{gtart_h}

\def\ifplaintex{\expandafter\ifx\csname documentclass\endcsname\relax}

\def\gtp{{\mathsurround=0pt\it $\cal G\mskip-2mu$eometry \&\ 
$\cal T\!\!$opology $\cal P\!$ublications}}  

\def\recd{{\small Received:\qua\receiveddate\ifx\reviseddate\relax
\else\qquad Revised:\qua\reviseddate\fi\par}} 

\def\lognumber#1{\def\thelognumber{#1}}
\def\volumenumber#1{\def\thevolumenumber{#1}}
\def\volumeyear#1{\def\thevolumeyear{#1}}
\def\papernumber#1{\def\thepapernumber{#1}}
\def\pagenumbers#1#2{\def\startpage{#1}\def\finishpage{#2}}
\def\published#1{\def\publishdate{#1}}

\def\received#1{\def\receiveddate{#1}}

\def\accepted#1{\def\accepteddate{#1}}

\long\def\asciiabstract#1{\long\def\theasciiabstract{#1}}

\let\\\par\let\thelognumber\relax\let\thevolumenumber\relax
\let\thepapernumber\relax\let\thevolumeyear\relax\let\startpage\relax
\let\finishpage\relax\let\publishdate\relax\let\receiveddate\relax
\let\reviseddate\relax\let\accepteddate\relax\let\theasciititle\relax
\let\theasciiauthors\relax
\let\theasciiabstract\relax

\let\theasciiemail\relax

\ifplaintex
\font\logobig=cmssbx10 scaled 3836
\font\logomed=cmssbx10 scaled 2557
\else
\font\logobig=cmssbx10 scaled 4200
\font\logomed=cmssbx10 scaled 2800
\fi

\long\def\makeagttitle{   
\count0=\startpage
\agt\hfill      
\hbox to 45truept{\vbox to 0pt{\vglue -13truept{\logomed A\kern -.37em{\logobig 
T}\kern -.38em G}\vss}\hss}
\break
{\small Volume \thevolumenumber\ (\thevolumeyear)
\startpage--\finishpage\nl
Published: \publishdate}

\vglue .25truein

{\parskip=0pt\leftskip 0pt plus
1fil\def\\{\par\smallskip}{\Large\bf\thetitle}\par\medskip} \vglue
0.05truein

{\parskip=0pt\leftskip 0pt plus 1fil\def\\{\par}{\sc\theauthors}
\par\medskip}%
 
\vglue 0.03truein 

{\small\leftskip 25truept\rightskip 25truept{\bf Abstract}\stdspace\theabstract

{\bf AMS Classification}\stdspace\theprimaryclass
\ifx\thesecondaryclass\relax\else; \thesecondaryclass\fi\par
{\bf Keywords}\stdspace \thekeywords\par}\vglue 7truept

}   

\ifplaintex
\font\phead=cmsl9 scaled 950
\font\pnum=cmbx10 scaled 913
\font\pfoot=cmsl9 scaled 950
\headline{\vbox to 0pt{\vskip -4.5mm\line{\small\phead\ifnum
\count0=\startpage ISSN 1472-2739 (on-line) 1472-2747 (printed)
\hfill {\pnum\folio}\else\ifodd\count0\def\\{ }%
\ifx\theshorttitle\relax\thetitle\else\theshorttitle\fi\hfill{\pnum\folio}
\else\def\\{ and }{\pnum\folio}\hfill\ifx\theshortauthors\relax\theauthors
\else\theshortauthors\fi\fi\fi}\vss}}
\footline{\vbox to 0pt{\vglue 0mm\line{\small\pfoot\ifnum\count0=\startpage
\copyright\ \gtp\hfill\else
\agt, Volume \thevolumenumber\ (\thevolumeyear)\hfill\fi}\vss}}
\else
\font=cmsl9 scaled 1050
\font\lnum=cmbx10 
\font=cmsl9 scaled 1050
\makeatletter
\def\@oddhead{{\small\lhead\ifnum\count0=\startpage ISSN 1472-2739 
(on-line) 1472-2747 (printed)\hfill {\lnum\number\count0}\else\ifodd\count0
\def\\{ }\ifx\theshorttitle\relax \thetitle \else\theshorttitle\fi\hfill
{\lnum\number\count0}\else\def\\{ and }{\lnum\number\count0}
\hfill\ifx\theshortauthors\relax 
\theauthors\else\theshortauthors\fi\fi\fi}}\def\@evenhead{\@oddhead}
\def\@oddfoot{\small\lfoot\ifnum\count0=\startpage\copyright\ \gtp\hfill\else
\agt, Volume \thevolumenumber\ (\thevolumeyear)\hfill\fi}
\def\@evenfoot{\@oddfoot}
\makeatother
\fi
\let\maketitlepage\makeagttitle
\let\makeshorttitle\maketitlepage
\let\maketitle\maketitlepage

\newwrite\gtoutfile
\long\gdef\makeheadfile{  
{\def\\{, }\def\s{ }
\immediate\openout\gtoutfile head.xxx
\immediate\write\gtoutfile{Proxy-for: \ifx\theasciiauthors\relax
\theauthors\else\theasciiauthors\fi\s<\ifx\theasciiemail\relax\theemail\else\theasciiemail\fi>}
\immediate\write\gtoutfile{\noexpand\\}
\immediate\write\gtoutfile{Authors: \ifx\theasciiauthors\relax
\theauthors\else\theasciiauthors\fi}
{\def\\{ }\immediate\write\gtoutfile{Title: \ifx\theasciititle\relax
\thetitle\else\theasciititle\fi}}
\immediate\write\gtoutfile{Subj-class: GT or SG, GR etc}
\immediate\write\gtoutfile{MSC-class: \theprimaryclass\ifx\thesecondaryclass\relax\else, \thesecondaryclass\fi}
\immediate\write\gtoutfile{Journal-ref: Algebr. Geom. Topol. \thevolumenumber\s
(\thevolumeyear) \startpage-\finishpage}
\immediate\write\gtoutfile{Comments: Published by Algebraic and
Geometric Topology at}
\immediate\write\gtoutfile{\s\s\s  http://www.maths.warwick.ac.uk/agt/AGTVol\thevolumenumber/agt-\thevolumenumber-\thepapernumber.abs.html}
\immediate\write\gtoutfile{\noexpand\\}
\immediate\write\gtoutfile{}
\ifx\theasciiabstract\relax
\immediate\write\gtoutfile{\theabstract}\else
\immediate\write\gtoutfile{\theasciiabstract}\fi
\immediate\write\gtoutfile{}
\immediate\write\gtoutfile{\noexpand\\}
\immediate\write\gtoutfile{}
\immediate\closeout\gtoutfile}}  

\def\maketitlepage{\makeagttitle\makeheadfile}
\let\makeshorttitle\maketitlepage
\let\maketitle\maketitlepage

\lognumber{54}
\volumenumber{5}
\volumeyear{2005}
\papernumber{54}
\pagenumbers{1389}{1418}
\received{15 July 2004} 
\accepted{8 July 2005}
\published{15 October 2005}

\usepackage{amsmath,amscd,amssymb,graphicx,psfrag}

\def\psfraga <#1,#2> #3#4{%
\psfrag {#3}{\smash{\rlap{\kern #1 \raise #2\hbox{#4}}}}}

\def\figref#1{\hyperlink{#1anchor}{Figure~\ref*{#1}}}
\def\anchor#1{\noindent\hypertarget{#1anchor}{\smash{$\phantom{99}$}}\newline}

\newtheorem{thm}{Theorem}[section]
\newtheorem{prop}[thm]{Proposition}
\newtheorem{cor}[thm]{Corollary}

\newtheorem{lem}[thm]{Lemma}

\theoremstyle{definition}

\newtheorem{remark}[thm]{Remark}

\newcommand{\ut}{\underline{t}}
\newcommand{\us}{\underline{s}}
\newcommand{\ui}{\underline{i}}

\newcommand{\lra}{\longrightarrow}

\newcommand{\Sig}{\Sigma}

\newcommand\Dual{\mathcal D}
\newcommand\Duality\Dual

\newcommand\relspinc{\underline{\spinc}}

\newcommand\x{\mathbf x}
\newcommand\w{\mathbf w}
\newcommand\z{\mathbf z}

\newcommand\y{\mathbf y}

\newcommand\ModSphere{\ModFlow\left({\mathbb S}\longrightarrow
\Sym^{g-1}(\Sigma_{1})\times \Sym^2(\Sigma_{2})\right)}
\newcommand\ModSpheres\ModSphere

\newcommand\Mas{\mu}
\newcommand\UnparModSp{\widehat \ModSp}
\newcommand\UnparModFlow\UnparModSp

\newcommand\PD{\mathrm{PD}}

\newcommand{\spinc}{\mathfrak s}

\newcommand{\sD}{\mathcal{D}}
\newcommand\CFLhat{\widehat{\CFL}}
\newcommand\sM{\mathcal{M}}
\newcommand\sF{\mathcal{F}}

\newcommand\HFLhat{\widehat{\HFL}}

\newcommand\ModMaps{\mathcal M}
\newcommand\ModSp\ModMaps

\newcommand\Ta{{\mathbb T}_{\alpha}}
\newcommand\Tb{{\mathbb T}_{\beta}}
\newcommand\Tc{{\mathbb T}_{\gamma}}

\newcommand\alphas{\mbox{\boldmath$\alpha$}}

\newcommand\betas{\mbox{\boldmath$\beta$}}
\newcommand\gammas{\mbox{\boldmath$\gamma$}}
\newcommand\thetas{\mbox{\boldmath$\theta$}}
\newcommand\deltas{\mbox{\boldmath$\delta$}}

\newcommand\spincrel\relspinc

\newcommand\HFK{\text{\it HFK\/}}
\newcommand\HFL{\text{\it HFL\/}}
\newcommand\CFL{\text{\it CFL\/}}
\newcommand\CFK{\text{\it CFK\/}}

\newcommand\BasePt{w}
\newcommand\FiltPt{z}

\newcommand{\R}{\mathbb{R}}

\newcommand{\Z}{\mathbb{Z}}

\newcommand{\ModSWfour}{\mathcal{M}}
\newcommand{\ModFlow}{\ModSWfour}

\newcommand{\SpinC}{{\mathrm{Spin}}^c}

\newcommand\abuts\Rightarrow
\newcommand\Sym{\mathrm{Sym}}

\title{Longitude Floer homology and the\\Whitehead double}
\author{Eaman Eftekhary}
\address{Mathematics Department, Harvard University\\1 Oxford Street, Cambridge, MA 02138, USA}
\email{eaman@math.harvard.edu}

\begin{document}

\begin{abstract}
We define the \emph{longitude Floer homology} of a knot $K\subset
S^3$ and show that it is a topological invariant of $K$. Some
basic properties of these homology groups are derived. In
particular, we show that they distinguish the genus of $K$. We
also make explicit computations for  the $(2,2n+1)$ torus knots.
Finally a correspondence between the longitude Floer homology of
$K$ and the Ozsv\'ath-Szab\'o Floer homology of its Whitehead
double $K_L$ is obtained.
\end{abstract}

\asciiabstract{%
We define the longitude Floer homology of a knot K in
S^3 and show that it is a topological invariant of K.  Some
basic properties of these homology groups are derived.  In
particular, we show that they distinguish the genus of K. We
also make explicit computations for  the (2,2n+1) torus knots.
Finally a correspondence between the longitude Floer homology of
K and the Ozsvath-Szabo Floer homology of its Whitehead
double K_L is obtained.}

\maketitle

\section{Introduction and main results}
Associated with a knot $K\subset S^3$, Ozsv\'ath and Szab\'o have
defined (\cite{OS-knot}) a series of of homology groups
$$\widehat{\HFK}(K), \HFK^\infty(K),\ \text{and} \  \HFK^\pm(K),$$
which are graded by $\SpinC$-structures $$\spinc \in
\SpinC(S^3_0(K))$$
 of the three-manifold $S^3_0(K)$, obtained by a zero surgery on
 $K$.

In this paper, first we introduce a parallel construction called the
\emph{longitude Floer homology}.

The construction of Ozsv\'ath and Szab\'o
 relies on, first finding a Heegaard diagram $$(\Sigma,\alpha_1,\ldots,\alpha_g;
\beta_2,\ldots,\beta_g)=(\Sigma,\alphas,\betas_0)$$ for the knot complement in $S^3$.
Then the meridian $\mu$ is
added as a special curve to obtain a Heegaard diagram $(\Sigma,\alphas,\{\mu\}\cup \betas_0)$
for $S^3$. A point $v$ on $\mu -\cup_i \alpha_i$
is specified as well. One will then put two points $z,w$ on the two sides of the curve
$\mu$ close to $v$. The  homology groups are  constructed as the Floer homology
associated with the totally real tori $\mathbb{T}_\alpha=\alpha_1 \times \ldots\times \alpha_g$
and $\mathbb{T}_\beta=\mu \times \beta_2 \times \ldots\times \beta_g$
in the symplectic manifold
$\Sym^g(\Sigma)$.

The marked points $z,w$ will filter the boundary maps and this leads to the Floer homology
groups $\widehat{\HFK}(K),\HFK^\infty(K)$ and $\HFK^\pm(K)$.

In this paper, instead of completing
$(\Sigma;\alpha_1,\ldots,\alpha_g;\beta_2,\ldots,\beta_g)$
to a Heegaard diagram for $S^3$ by adding the meridian of the knot $K$, we choose the special
curve $\hat \beta_1$ to be the longitude of  $K$ sitting on the surface $\Sigma$
and not cutting any of the $\beta$ curves. This choice is made so that it
leads to a Heegaard diagram
$(\Sigma;\alpha_1,\ldots,\alpha_g;\hat \beta_1,\beta_2,\ldots,\beta_g)$ of the three-manifold
$S^3_0(K)$, obtained by a zero surgery on the knot $K$. Similar to the construction of
Ozsv\'ath and Szab\'o we also fix a marked point $v$ on the curve $\hat \beta_1 -\cup_i
\alpha_i$ and choose the base points $z,w$ next to it, on the two sides of $\hat \beta_1$.

We call the resulting Floer homology groups the \emph{longitude Floer homologies} of the
knot $K$, and  denote them by $\widehat{\HFL}(K),\HFL^\infty(K)$ and $\HFL^\pm(K)$. These are
graded by a $\SpinC$ grading $\underline{s}\in \tfrac{1}{2}+\mathbb{Z}$ and a (relative)
Maslov grading $m$. Generally, the $\SpinC$ classes in
$\SpinC(S^3_0(K))$ will be denoted by $\spinc, \spinc_1$, etc, while
if an identification $\Z \simeq \SpinC(S^3_0(K))$ is fixed these
classes are denoted by $s,s_1$, etc.  Similarly, in general classes in
$\tfrac{1}{2}\text{PD}[\mu]+\SpinC(S^3_0(K))$ are denoted by
$\relspinc,\relspinc_1$, etc, while under a fixed identification with
$\tfrac{1}{2}+\Z$ they are denoted by $\underline{s},\underline{s}_1$,
etc.

We show that the following holds for this theory (cf.\ \cite{OS-genus}):

\begin{thm}
{ Suppose that $g$ is the genus of the nontrivial knot $K\subset S^3$. Then
\begin{displaymath}
\widehat{\HFL}(K,g-\tfrac{1}{2})\simeq \widehat{\HFL}(K,-g+\tfrac{1}{2})\neq 0,
\end{displaymath}
and for any $\ut>g-\tfrac{1}{2}$ in $\tfrac{1}{2}+\mathbb{Z}$ the groups
$\widehat{\HFL}(K,\ut)$ and $\widehat{\HFL}(K,-\ut)$ are both trivial. Furthermore,
for any $\ut$ in $\tfrac{1}{2}+\mathbb{Z}$ there is an isomorphism of the
relatively graded groups (graded by the Maslov index)
\begin{displaymath}
\widehat{\HFL}(K,\ut)\simeq \widehat{\HFL}(K,-\ut).
\end{displaymath}
}
\end{thm}

Using some other results of Ozsv\'ath and Szab\'o we also show that:

\begin{cor}
{ If $K$ is a fibered knot of genus $g$ then
$\widehat{\HFL}(K,\pm(g-\frac{1}{2}))$ will be equal to
$\mathbb{Z}\oplus{\mathbb{Z}}$.
}\end{cor}

In the second part of this paper  we will construct a Heegaard diagram
of the Whitehead double
of a knot $K$ and will show a correspondence between the
longitude Floer homology of the given knot $K$ and the Ozsv\'ath-Szab\'o Floer homology of
the Whitehead double of $K$ in the non-trivial $\SpinC$ structure.

The Whitehead double of $K$ is a special case of a construction called the
\emph{satellite construction}. Suppose that
$$
e:D^2\times S^1 \rightarrow S^3$$
is an embedding such that $e(\{0\}\times S^1)$ represents the knot $K$,
and such that $e(\{1\}\times S^1)$ has zero linking number with
$e(\{0\}\times S^1)$. Let
$L$ be a knot in $D^2\times S^1$, and denote $e(L)$ by $K_L$. $K_L$ is called
a satellite of the knot $K$. 

The Alexander polynomial of $K_L$ may be easily expressed
in terms of the Alexander polynomial of $K$ and $L$. In fact, if $L$ represents $n$ times
the generator of the first homology of $D^2\times S^1$, and $\Delta_K(t)$ and
$\Delta_L(t)$ are the Alexander polynomials of $K$ and $L$, then the Alexander polynomial
of $K_L$ is given by (see [1] for a proof)
\begin{equation}
\Delta_{K_L}(t)=\Delta_K(t^n).\Delta_L(t).
\end{equation}
\begin{figure}[ht!]\small\anchor{fig:whitehead.1}
\psfrag{L}{$L$}
\cl{\includegraphics[width=2.2in]{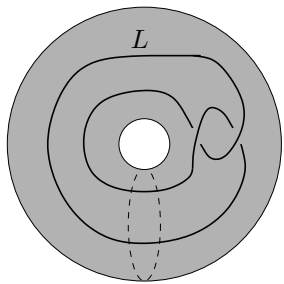}\hspace{3cm}}
\caption{\label{fig:whitehead.1}
The knot $L$ is used in the satellite construction to obtain the
Whitehead double of a knot $K$}
\end{figure}
In particular, if $L$ is an embedding of the
un-knot in $D^2\times S^1$ which is shown in  \figref{fig:whitehead.1},
$K_L$ is called the Whitehead double of $K$.

The Alexander polynomial of the Whitehead double of a knot $K$
is always trivial for this choice  of the
framing for a knot $K$, given by $e(\{1\} \times S^1)$.
It is known  that the Whitehead double of a nontrivial knot,
is nontrivial (see \cite{Burde}, for a proof). So the genus of $K_L$ is at least $1$. There is a surface of genus $1$
in $D^2\times S^1$ which bounds the knot $L$. The image of this surface under the map $e$
will be a Seifert surface of genus $1$ for $K_L$. This shows that $g(K_L)=1$.
By the result of \cite{OS-genus}, $\widehat{\HFK}(K_L,\pm 1)$ is nontrivial and
 $\widehat{\HFK}(K_L,s)=0$ for $|s|>1$.

We will show that the Ozsv\'ath-Szab\'o Floer homology (the hat theory)
in the $\SpinC$ structures $\pm 1$ is
in fact closely related to the longitude Floer homology introduced in
the first part.

More precisely we show:

\begin{thm}
{Let $K_L$ denote the Whitehead double of a knot $K$ in $S^3$. The
Ozsv\'ath-Szab\'o Floer homology group $\widehat{\HFK}(K_L,1)$ is
isomorphic to
 $$\widehat{\HFL}(K)=\bigoplus_{\ui\in \mathbb{Z}+\frac{1}{2}}\widehat{\HFL}(K,\ui)$$ as (relatively)
$\mathbb{Z}$-graded abelian groups with the (relative) grading on both sides coming from the
Maslov grading.}
\end{thm}

\textbf{Acknowledgment}\qua This paper is part of my PhD thesis in
Princeton, and I would like to use this opportunity to thank my
advisor Zolt\'an Szab\'o for all his help and support. I would also
like to thank Peter Ozsv\'ath for helpful discussions, Jake Rasmussen
for pointing out a mistake in an earlier version of this paper, and
the referee for useful comments.

\section{Construction of longitude homology}

Suppose that an oriented  knot $K$ is given in $S^3$. We may consider a Heegaard diagram
$(\Sigma_g;\alpha_1,\ldots,\alpha_g;\hat \beta_1,\beta_2,\ldots,\beta_g;v)$ such that
$\Sigma_g$ is a Riemann surface of genus $g$ and $\alphas=(\alpha_1,\ldots,\alpha_g)$
and $\betas=(\hat \beta_1,\beta_2,\ldots.,\beta_g)=\{\hat \beta_1\}\cup \betas_0$ are two $g$-tuples of disjoint
simple closed
loops on $\Sigma_g$ such that the elements of each $g$-tuple are linearly
independent in the first homology of $\Sigma_g$. Here $v$ is a marked point
on $\hat \beta_1 -\alphas$.
We assume furthermore that
$(\Sigma_g;\alpha_1,\ldots,\alpha_g;\beta_2,\ldots,\beta_g)$ is a Heegaard diagram
for the complement of the knot $K$ in $S^3$, and that $\hat \beta_1$ represents
the oriented longitude of the knot $K$. We assume that $\hat \beta_1$ is chosen so that
 the whole Heegaard diagram
$$(\Sigma_g;\alpha_1,\ldots,\alpha_g;\hat \beta_1,\beta_2\ldots,\beta_g)=(\Sigma_g;\alphas;\betas)$$
is a Heegaard diagram for the three-manifold $S^3_0(K)$ obtained by a zero surgery  on the knot $K$.

Choose two base points $z,w$ in the complement
$$\Sigma_g - \alphas-\betas,$$
very close to the marked point $v$ on $l=\hat \beta_1$, such that $z$ is on  the right
hand side and $w$ is on the left hand side of $l$ with respect to the
orientation of $l$ and that of $\Sigma$ coming from the handlebody determined
by $(\mu,\beta_2,\ldots,\beta_g)$. Here $\mu$ represents the meridian of $K$.
The usual construction of
Ozsv\'ath and Szab\'o works to give us a well defined Floer homology theory
associated with this setup.

Namely we may consider the two tori
$$\mathbb{T}_\alpha=\alpha_1\times \ldots.\times \alpha_g,
\mathbb{T}_\beta=\hat \beta_1 \times \ldots\times \beta_g \subset \Sym^g (\Sigma_g)$$
as two totally real subspaces of the symplectic manifold $\Sym^g (\Sigma_g)$, which
is the $g$-th symmetric product of the surface $\Sigma_g$. If the curves are transverse
on $\Sigma_g$ then these two $g$ dimensional submanifolds will intersect each other
transversely in finitely many points. The complexes $\CFL^\infty, \CFL^{\pm}$ and
$\widehat{\CFL}$ are generated
by the generators $[\x,i,j], i,j\in \mathbb{Z}$, in the infinity
theory, by $[\x,i,j], i,j\in \mathbb{Z}^{\leq 0}$ in $\CFL^-$ case and by the elements
$[\x,0,0]$ in the hat theory, where $\x$ is an intersection point of
$\mathbb{T}_\alpha$ and $ \mathbb{T}_\beta$. i.e.\
$\x\in \mathbb{T}_\alpha \cap \mathbb{T}_\beta$. The groups $\CFL^+$ will appear as the
quotient of the embedding
$$
0\longrightarrow \CFL^-(K) \longrightarrow \CFL^\infty(K).$$
For any two intersection points $\x,\y \in \mathbb{T}_\alpha \cap \mathbb{T}_\beta$
there is a well-defined element $\epsilon (\x,\y) \in H_1(S^3_0(K))$ defined as follows:

Choose a path $\gamma_\alpha$ from $\x$ to $\y$ in $\mathbb{T}_\alpha$ and a path
$\gamma_\beta$ from $\y$ to $\x$ in $\mathbb{T}_\beta$. $\gamma_\alpha +\gamma_\beta$
will represent an element of $H_1(\Sym^g(\Sigma_g),\mathbb{Z})$ which is well-defined
modulo $H_1(\mathbb{T}_\alpha, \mathbb{Z})\oplus H_1(\mathbb{T}_\beta,\mathbb{Z})$.
Thus, we will get our desired element
\begin{displaymath}
\epsilon(\x,\y)=[\gamma_\alpha+\gamma_\beta]\in
\frac{H_1(\Sym^g(\Sigma_g),\mathbb{Z})}
{H_1(\mathbb{T}_\alpha, \mathbb{Z})\oplus H_1(\mathbb{T}_\beta,\mathbb{Z})}
\cong H_1(S^3_0(K),\mathbb{Z}).
\end{displaymath}
This first homology group is $H_1(S^3_0(K),\mathbb{Z})\cong \mathbb{Z}$, so we will
get a relative $\mathbb{Z}$-grading on the set of generators. There will be  maps
$$\spinc_z,\spinc_w:\mathbb{T}_\alpha \cap \mathbb{T}_\beta \rightarrow \SpinC(S^3_0(K))$$
as in \cite{OS-knot,OS-3m1},
which may be defined using each of the points $z$ or $w$. Note that our maps $\spinc_z,\spinc_w$
were called $\relspinc_z,\relspinc_w$ in \cite{OS-knot,OS-3m1}.
Unlike the standard Heegaard Floer homology of Ozsv\'ath and Szab\'o, $\spinc_z(\x)$ does not
agree with the $\SpinC$ structure $\spinc_w(\x)$. However, we will have
$$\spinc_z(\x)=\spinc_w(\x)+\text{PD}[\mu],$$ where $\mu $ is the meridian of the knot $K$,
thought of as a closed loop in $S^3_0(K)$, generating its first
homology (cf.\ \cite{OS-knot}, lemma 2.19). We may either decide
to assign the $\SpinC$ structure $\spinc_z(\x)$ to $\x$, or more invariantly define
$$\relspinc(\x)=\frac{\spinc_z(\x)+\spinc_w(\x)}{2}\in \tfrac{1}{2}\text{PD}[\mu]+\SpinC(S^3_0(K))\simeq
\tfrac{1}{2}+\Z.$$
We will frequently switch between these two choices, distinguishing them by the name of the maps, i.e.\
$\spinc_z$ versus $\relspinc$.
Also, the elements of $\SpinC(S^3_0(K))$ will be denoted by $\spinc,\spinc_1$, etc,
while the elements of $\frac{1}{2}\PD[\mu]+\SpinC(S^3_0(K))$ will be denoted by
$\relspinc,\relspinc_1$, etc (and if an identification with $\Z$
(resp. $\tfrac{1}{2}+\Z$) is fixed, by $s,s_1$,
etc. (resp. $\underline{s},\underline{s}_1$, etc)).

If $\epsilon(\x,\y)=0$, which is the same as $\spinc_z(\x)=\spinc_z(\y)$, then there is a family of
homotopy classes of disks with boundary in $\mathbb{T}_\alpha$ and $ \mathbb{T}_\beta$,
connecting  $\x$ and $\y$. Note that to each map
\begin{displaymath}
\begin{split}
&u: [0,1]\times \mathbb{R} \rightarrow \Sym^g(\Sigma_g)\\
&u(\{0\} \times \mathbb{R}) \subset \mathbb{T}_\alpha, \ \ \
u(\{1\} \times \mathbb{R}) \subset \mathbb{T}_\beta\\
&u(s,t) \rightarrow \x \ \  \text{as} \ \ t\rightarrow \infty, \ \ \ \ \
u(s,t) \rightarrow \y \ \  \text{as} \ \ t\rightarrow -\infty,
\end{split}
\end{displaymath}
we may assign a domain $\mathcal{D}(u)$ on $\Sigma_g$ as follows: Choose a point $z_i$ in
$D_i$, where $D_i$'s are the connected components of the complement of the curves:
$$\Sigma_g-\alphas-\betas=\Sigma_g-\alpha_1-\ldots-\alpha_g-\hat \beta_1-\ldots-\beta_g=\coprod_i D_i. $$
Let $L_i$ be the codimension $2$ subspace $\{z_i\}\times \Sym^{g-1}(\Sigma_g)$ of
$\Sym^g(\Sigma_g)$. We denote the intersection number of $u$ and $L_i$ by $n_{z_i}(u)$.
This number only depends on the homotopy type of the disk $u$, which we denote
by $\phi$. We will get a formal sum of domains
$$
\mathcal{D}(\phi)=\mathcal{D}(u)=\sum_{i=1}^{n} n_{z_i}(u).D_i.$$
For any two points $\x,\y\in \mathbb{T}_\alpha \cap \mathbb{T}_\beta$ with $\epsilon(\x,\y)=0$,
there is a unique homotopy type $\phi$ of the disks connecting
them with $n_z(\phi)=i$ and $n_w(\phi)=j$.
These are described
as follows. If $\phi, \psi$ denote two different homotopy types of disks between
$\x$ and $\y$, then $\mathcal{D}(\phi) -\mathcal{D}(\psi)$ will be a domain whose boundary is
a sum of the curves $\alpha_1,\ldots,\alpha_g,\hat \beta_1,\beta_2,\ldots,\beta_g$.
These are called \emph{periodic domains}. Note that Ozsv\'ath and
Szab\'o require that a periodic domain should have zero multiplicity
at one of the prescribed marked points (cf.\ \cite{OS-3m1}), while we allow
any arbitrary multiplicities at the marked points $z$ and $w$. 

The set of periodic domains is isomorphic to the torsion free part of
$H_1(Y,\mathbb{Z})$ plus $\mathbb{Z}$
(where $Y$ is the three manifold determined by the Heegaard
diagram). This is $H_1(S^3_0(K),\mathbb{Z})\oplus \mathbb{Z}$
in the above case. As a result, the set of periodic domains
is isomorphic to $\mathbb{Z}\oplus \mathbb{Z}$ for this problem. This way, we may
assign an element in
$H_1(S^3_0(K),\mathbb{Z})\oplus \mathbb{Z}\cong \mathbb{Z}\oplus \mathbb{Z}$ to any pair $\phi,\psi$ of homotopy
disks between $x,y$ which will be denoted by $h(\phi,\psi)$.

Denote the generators of the set of periodic domains by $\mathcal{D},\mathcal{D}_0$.
Here $\mathcal{D}_0$ is the disk whose domain is the whole surface
 $\Sigma$ and $\mathcal{D}$ is characterized with the property that
$n_z(\mathcal{D})=0$ and $ n_w(\mathcal{D})=1$. 

Associated with each homotopy class $\phi$ of disks connecting $\x$ and $\y$, which is
denoted by $\phi \in \pi_2(\x,\y)$, we define
\begin{displaymath}
\mathcal{M}(\phi)=\Big
\{\begin{array}{ccc}
u:[0,1]\times \mathbb{R} \rightarrow \Sym^g(\Sigma_g) &
\Big |&
u\in \phi, \ \ \ \ \bar \partial_{J_t} u(s,t)=0
\end{array}\Big \},
\end{displaymath}
where $J=\{J_t\}_{t\in [0,1]}$ is a generic one parameter family of almost complex structures
 arising from  complex structures on the
surface $\Sigma_g$ (cf.\ \cite{OS-3m1}). This moduli space will have an expected dimension which
we will denote by $\mu(\phi)$ ($\mu(\phi)$ should not be confused with the meridian
$\mu$ of the knot $K$).

There is a $\mathbb{R}$-action, as usual, and we may divide the moduli space by
this action. Denote the quotient $\mathcal{M}(\phi)/\mathbb{R}$ of this moduli space
by $\UnparModFlow(\phi)$.
The boundary maps for the infinity theory $\CFL^\infty (K)$ are
described as follows. Let $\CFL^\infty(K,\relspinc)$ be the part of the complex $\CFL^\infty(K)$
generated by the intersection points $\x$ such that $\relspinc(\x)=\relspinc \in \frac{\text{PD}[\mu]}{2}+\SpinC(S^3_0(K))$. Define
the boundary map $\partial^\infty$ for a generator $[\x,i,j]$,
$\x\in \mathbb{T}_\alpha \cap \mathbb{T}_\beta$ with $\relspinc(\x)=\relspinc$, to be the sum
\begin{displaymath}
\partial^\infty[\x,i,j]=\sum_{\substack{\y\in\Ta\cap\Tb}}
                              \sum_{\substack{\phi\in\pi_2(\x,\y)\\
                    \Mas(\phi)=1}}
\#\left(\UnparModFlow(\phi)\right)[\y,i-n_{\BasePt}(\phi),j-n_{\FiltPt}(\phi)].
\end{displaymath}
If the diagram is strongly admissible (cf.\ \cite{OS-3m1,OS-3m2}), then the above sum is always finite.

Let $\partial^{-}$ be the restriction to the complex $\CFL^{-}(K,\relspinc)$ and $\partial^+$ to be
the induced map on the quotient complex $\CFL^+(K,\relspinc)$. On $\widehat{\CFL}(K,\relspinc)$ we consider
the simpler map:
\begin{displaymath}
\partial[\x]=\sum_{\substack{\y\in\Ta\cap\Tb}}
             \sum_{\substack{ \phi\in\pi_2(\x,\y)\\
                 \Mas(\phi)=1\\
                  n_z(\phi)=n_w(z)=0}}
\#\left(\UnparModFlow(\phi)\right)[\y].
\end{displaymath}
Here there is no need for an admissible Heegaard diagram, since there will be
finitely many terms involved in the above sum.

These maps are all differentials which compose with themselves to give zero
(the proof is identical with those used by Ozsv\'ath and Szab\'o).
As a result we will get the homology groups
$$
\HFL^\infty(K,\relspinc)\ \  ,\ \  \HFL^\pm(K,\relspinc)\ \ ,\text{ and }
\ \widehat{\HFL}(K,\relspinc).$$
Here $\relspinc$ is taken to be in
$$\tfrac{1}{2}\text{PD}[\mu]+\SpinC(S^3_0(K))\simeq  \tfrac{1}{2}+\Z.$$
We should prove that these homology groups are independent of the choice
of the Heegaard diagram and the almost complex structure and that they only
depend on the knot $K$ and the $\SpinC$-structure $\relspinc$ chosen from
$\frac{1}{2}\text{PD}[\mu]+\SpinC(S^3_0(K))$. The independence from the almost complex structure
is proved exactly in the same way that the knot Heegaard Floer homology
is proved to be independent of this choice.

\begin{thm}
Let $K$ be an oriented knot in $S^3$ and fix the
$\SpinC$-structure $\relspinc\in \frac{1}{2}\text{PD}[\mu]+\SpinC(S^3_0(K))$. Then the homology groups
$\HFL^\infty(K,\relspinc),\HFL^\pm(K,\relspinc)$, and $\widehat{\HFL}(K,\relspinc)$ will be topological invariants
of the oriented knot $K$ and  the $\SpinC$-structure $\relspinc$;
i.e.\ They are independent of the choice of the marked Heegaard diagram $$(\Sigma_g;
\alpha_1,\ldots,\alpha_g;\hat \beta_1,\beta_2,\ldots,\beta_g;v)$$ used in the definition.
\end{thm}
\begin{proof} The proof is almost identical with the proof in the case of Heegaard Floer homology.
We will just sketch the steps of this proof.
We remind the reader of the following proposition (prop. 3.5.) of \cite{OS-knot}:

\begin{prop}{If two Heegaard diagrams
$$
(\Sigma;\alpha_1,\ldots,\alpha_g;\hat \beta_1,\beta_2,\ldots,\beta_g),
(\Sigma;\alpha_1',\ldots,\alpha_g';\hat \beta_1',\beta_2',\ldots,\beta_g'),$$
represent the same manifold obtained by zero surgery on the knot $K$ in $S^3$, then
we can pass from one to the other by a sequence of the following  moves and
their inverses:

{\rm(1)}\qua Handle slide and isotopies among 
$\alpha_1,\ldots,\alpha_g,\beta_2,\ldots,\beta_g$.

{\rm(2)}\qua Isotopies of $\hat \beta_1$.

{\rm(3)}\qua Handle slides of $\hat \beta_1$ across some of the $\beta_2,\ldots,\beta_g$.

{\rm(4)}\qua Stabilization (introducing cancelling pairs
$\alpha_{g+1},\beta_{g+1}$ and increasing the genus of $\Sigma$ by
one).}
\end{prop}

As in \cite{OS-knot} assume that
$$
D_1=(\Sigma;\alpha_1,\ldots,\alpha_g;\beta_1,\ldots,\beta_g;z,w),\ \ \
D_2=(\Sigma;\beta_1,\ldots,\beta_g;\gamma_1,\ldots,\gamma_g;z,w)$$
are a pair of doubly pointed Heegaard diagrams. There will be a map:
$$F:\CFL^\infty (D_1)\otimes \CFL^\infty(D_2) \longrightarrow
\CFL^\infty(\Sigma;\alpha_1,\ldots,\alpha_g;\gamma_1,\ldots,\gamma_g;z,w)
$$
defined by
\begin{equation}
\begin{split}
F(\partial[\x&,i,j] \otimes [\y,l,k])=\\
 &\sum_{\substack{\z \in\Ta\cap\Tc}}\sum_{\substack{
                  \phi\in\pi_2(\x,\y,\z)\\
          \Mas(\phi)=0}}
\#\left(\UnparModFlow(\phi)\right)[\z,i+l-n_z(\phi),j+k-n_w(\phi)].
\end{split}
\end{equation}
Here we use the notation $\pi_2(\x,\y,\z)$ for the space of homotopy classes of the
disks $u:\Delta \rightarrow \Sym^g(\Sigma_g)$ from the unit triangle
$\Delta$ with edges $e_1,e_2,e_3$, to the symmetric space, such that
$$
u(e_1)\subset \mathbb{T}_\alpha, \ \ \ u(e_2)\subset \mathbb{T}_\beta,
\ \ \ u(e_3)\subset \mathbb{T}_\gamma,$$
and the vertices of $\Delta$ are mapped to the three points $\x,\y,\z$.

Back to the proof of the theorem, the independence from the isotopies of
$\beta_2,\ldots,\beta_g$, the isotopies of $\alpha_1,\ldots,\alpha_g$
and even for the isotopies of the special curve $\hat \beta_1$ are easy
and identical to the standard case. The same
is true for handle slides among $\beta_2,\ldots,\beta_g$. In fact if
$\beta_2',\ldots,\beta_g'$ are obtained from $\beta_2,\ldots,\beta_g$ by a
handle slide and if $\hat \beta_1'$ is a small perturbation of
$\hat \beta_1$ then in the Heegaard diagram
$$
(\Sigma_g;\hat \beta_1,\beta_2,\ldots,\beta_g;
\hat \beta_1',\beta_2',\ldots,\beta_g';z,w),$$
$z,w$ lie in the same connected component of  complement of the
curves $$\Sigma-\hat \beta_1-\ldots-\beta_g-\hat \beta_1'-\ldots-\beta_g'.$$
We may assume that this is a strongly admissible Heegaard diagram
for $\#^g(S^2\times S^1)$ for the $\SpinC$-structure $\spinc_0$ on
$\#^g(S^2\times S^1)$ with trivial first Chern class. Denote by
$\CFL^\infty_\delta(\Sigma,\betas,\betas')$ the complex generated by the
generators $[\x,i,i]$, where $\betas$ represents
$(\hat \beta_1,\beta_2,\ldots,\beta_g)$ and $\betas'$ represents
$(\hat \beta'_1,\beta'_2,\ldots,\beta'_g)$. Then
\begin{displaymath}
\HFL^{\leq 0}_\delta(\Sigma;{\betas};{\betas'};z,w)\simeq
\mathbb{Z}[U]\otimes_\mathbb{Z} \Lambda^* H_1(T^g),
\end{displaymath}
where $U$ is the map sending $[\x,i,j]$ to $[\x,i-1,j-1]$. There will be
a top generator $\Theta \in \Lambda^*H_1(T^g)$ in the $\SpinC$ class
with trivial associated first Chern class. We  may use
\begin{displaymath}
[\x,i,j] \mapsto F([\x,i,j]\otimes \Theta)
\end{displaymath}
to define the map associated with a handle slide determined
by the pair $(\betas,\betas')$. This induces a map
in the level of homology which we may argue -- as is typical in the
previous work of Ozsv\'ath and Szab\'o (see \cite{OS-3m1}) -- that  in fact
induces an isomorphism.
The induced map on the subcomplex $\CFL^-(K)$, and also on the
complexes $\CFLhat^+(K)$ and $\CFLhat(K)$ will be isomorphisms as well,
since the map $F$ respects the filtration of $\CFL^\infty(K)$.
See \cite{OS-knot} for more
details. Other handle slides are quite similar. The proof of
independence from the handle addition is identical to the standard case.
\end{proof}

\begin{remark} This theory is in fact an extension of the usual Heegaard Floer homology for the
meridian of the knot $K\subset S^3$, considered as an image of $S^1$
 in $S^3_0(K)$. The meridian $\mu$ is not null-homologous
in $S^3_0(K)$ which makes it not satisfy the requirements of the construction of
\cite{OS-knot}.
However, the only difference that is forced to us, is that the maps $\spinc_z$ and
$\spinc_w$ will not assign the same $\SpinC$ structure to the generators of the complex.
As we have seen this is not a serious problem at all.
The independence of the homology groups from the choice of a
Heegaard diagram may be proved yet, as noted above.
\end{remark}

\section{Basic properties}

In this section we start developing some properties of these longitude Floer homology groups.
Let $K$
be a knot in $S^3$ and
$$
(\Sigma,\alpha_1,...,\alpha_g,\hat \beta_1,\beta_2,...,\beta_g;z,w)
$$
be a Heegaard diagram for $K$; where $z$ and $w$ are on the two sides of the longitude
$\hat \beta_1$. We may assume that the meridian $m$ of the knot is a curve on
$\Sigma$ which cuts $l=\hat \beta_1$ and one of the $\alpha$ curves, say $\alpha_1$,
exactly once, and is disjoint from all other curves $\alpha_i$ and $\beta_i, i\geq 2$.
 Denote the unique intersection point between $m$ and $\alpha_1$ by $x$. 

\begin{figure}[ht!]\small\anchor{fig:spin1}
\psfrag {1}{$l$}
\psfraga <-2pt,2pt> {m}{$m$}
\psfraga <-1pt,0pt> {z}{$z$}
\psfraga <-2pt,0pt> {w}{$w$}
\psfrag {a}{$\alpha_1$}
\psfraga <0pt,3pt> {x}{$x$}
\psfraga <-6pt,0pt> {x1}{$x_1$}
\psfraga <-3pt,0pt> {x2}{$x_2$}
\psfraga <-3pt,0pt> {x3}{$x_3$}
\psfrag {y1}{$y_2$}
\psfraga <-3pt,2pt> {y2}{$y_3$}
\psfrag {y3}{$y_1$}
\cl{\hspace{.8cm}\includegraphics{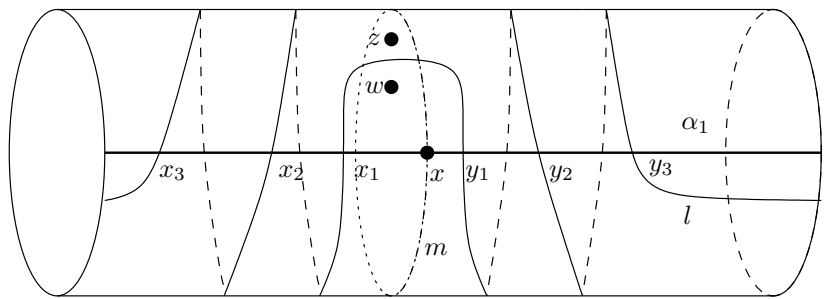}}
\caption{\label{fig:spin1}
{Let the curve $\beta_1$ wind around  the meridian curve $m$ sufficiently many times
and put $z$ and $w$ near $l$ in the inner most regions, and on the two sides of
$l$. }}
\end{figure}

Choose a large number $N$ and change the curve $l$ by winding it $N$
times around $m$ (cf.\ \cite{OS-3m1}, section 5). This will also be a Heegaard diagram for the same knot $K$ and we may assume
that the base points $z$ and $w$ are in the inner-most regions, as is 
shown in \figref{fig:spin1}.

There will be $2N$ new intersection points $x_1,...,x_N,y_1,...,y_N$ created   between
the two curves $\alpha_1,l$.

There is a periodic domain for the Heegaard diagram $(\Sigma,\alphas,\betas,z)$
which has multiplicity $1$ on the region containing $w$ and multiplicity zero at $z$.
The multiplicities of the
domains outside the cylinder shown in the figure will be  negative numbers less than
some fixed number $-N+k$. By choosing $N$ large enough, we may assume that this number
is sufficiently negative.

Remember that the generators of the complex $\widehat{\CFK}(K)$,
when computed using the standard
Heegaard diagram associated with $K$ (which comes from a knot projection), are in one-to-one
correspondence with combinatorial objects called the \emph{Kauffman states}
(see \cite{OS-knot,
OS-alternating} for more details). We will abuse the language and some times use the word
Kauffman state to refer to the generators of the chain complexes.

The Kauffman states of the above Heegaard diagram are of two types:

(1)\qua Those of the form $\{x_i,\bullet\}$ or $\{y_i,\bullet\}$, which are in correspondence
with the Kauffman state $\{x,\bullet\}$ of the Heegaard diagram
$$
(\Sigma,\alpha_1,...,\alpha_g;m,\beta_2,...,\beta_g;z)
$$
for the sphere $S^3$.

(2)\qua Those which are not of this form; We will  call them \emph{bad} Kauffman states.

There is a $\SpinC$ structure of $S^3_0(K)$ assigned to each Kauffman state using the base
point $z$.
Any two Kauffman states of the form  $\x=\{x_i,\bullet\}$ and
$\y=\{y_i,\bullet\}$ are in the same $\SpinC$-class
$\relspinc(\x)=\relspinc(\y)$. There is a difference of $\ell.[\Delta]$
between $\relspinc(\{x_i,\bullet\})$ and $\relspinc(\{x_{i+\ell},\bullet\})$, where $[\Delta]$
is the generator of the second homology group of $S^3_0(K)$.

We may choose the fixed number $k$ so that the $\SpinC$ difference between a bad Kauffman state
and a Kauffman state of the form $\{x_i,\bullet\}$ is at least $(N-i-k)[\Delta]$.

Let $\mathcal{D}$ be the periodic domain considered above, and let $\mathcal{D}_0$
denote the  periodic domain represented by the surface $\Sigma$. The space of all
periodic domains is generated by $\mathcal{D}$ and $\mathcal{D}_0$.

\begin{figure}[ht!]\small\anchor{fig:spin2}
\psfrag {1}{$1$}
\psfrag {2}{$2$}
\psfrag {3}{$3$}
\psfrag {x}{$x_3$}
\psfrag {y}{$y_3$}
\cl{\hspace{-.7cm}\includegraphics{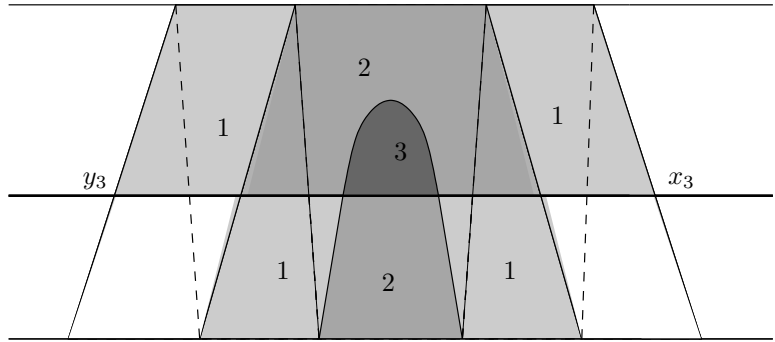}}
\caption{\label{fig:spin2}
{There is a domain connecting $\{x_i,\bullet\}$ and $\{y_i,\bullet\}$ with Maslov index
$1$ and coefficients $i-1,i$ at $z,w$ respectively. The domain for $\{x_3,\bullet\}$
and $\{y_3,\bullet\}$ is illustrated.}}
\end{figure}

There is a disk between $\{x_i,\bullet\}$ and $\{y_i,\bullet\}$ with Maslov index
$1$ and coefficients $i-1$ and $i$ at $w$ and $z$ respectively. This domain is illustrated for
$\{x_3,\bullet\}$ and $\{y_3,\bullet\}$ in  \figref{fig:spin2}. Let us denote this domain
by $\mathcal{D}_i$. Then
$$
\tilde{\mathcal{D}}_i=\mathcal{D}_i-\mathcal{D}-(i-1)\mathcal{D}_0
$$
 will be
the unique connecting domain between $\{x_i,\bullet\}$ and $\{y_i,\bullet\}$ with
zero coefficients on $z$ and $w$.
Note that the Maslov index of the domain $\mathcal{D}_0$ is equal to $2$. As a result,
\begin{equation}
\mu(\tilde{\mathcal{D}}_i)=1-2(i-1)-\mu(\mathcal{D}).
\end{equation}
$\mathcal{D}$ represents the generator of the second homology of $Y=S^3_0(K)$. Namely,
we may think of the Heegaard diagram for $Y$ as given by a Morse function $h$, and
assume that
\begin{displaymath}
\partial \mathcal{D}=\sum_{i=1}^{g}n_i\alpha_i+\sum_{i=1}^{g}m_i\beta_i.
\end{displaymath}
The points that flow to $\alpha_i$ form a disk $P_i$ that caps $\alpha_i$. Similarly, the points
that lie on the flow  coming out of $\beta_i$ form another disk $Q_i$ that caps the curve
$\beta_i$. Then with an appropriate orientation on $P_i$ and $Q_i$, so that $\partial P_i=-\alpha_i$
and $\partial Q_i=-\beta_i$, the domain
$$
\mathcal{D}+\sum_{i=1}^{g}n_i.P_i+\sum_{i=1}^{g}m_i.Q_i
$$
will represent a homology class $[F]$ in the three-manifold $Y=S^3_0(K)$ which in fact
generates its second homology (and so is equal to $\Delta$).

The Maslov index of this homology class is equal to
\begin{equation}
\chi({\mathcal{D}})=\langle c_1(s_i),[F]\rangle ,
\end{equation}
where $\chi(\mathcal{D})$ is the Euler measure of $\mathcal{D}$, and $s_i$ is the $\SpinC$
structure,
\begin{displaymath}
s_i=\spinc_z(\{x_i,\bullet\})=\spinc_z(\{y_i,\bullet\})\in \Z=\SpinC(S^3_0(K)).
\end{displaymath}
Again, this last identification is done so that the $\SpinC$ class
with trivial first Chern class is identified with $0\in \Z$.

As a result, $\mu(\tilde{\mathcal{D}}_i)=1-2(i-1)-\langle c_1(s_i),[F] \rangle$.
This domain has
very positive coefficients in the domains of the surface $\Sigma$, which are not on the
cylinder shown in the figure, if the index $i$ is not very big. The Maslov index
is $1$ exactly when
\begin{displaymath}
-2(i-1)=\langle c_1(s_i),[F]\rangle .
\end{displaymath}
Note that
$c_1(s_i)=-2(i-1)\text{PD}[ F]+c_1(s_1)$, which implies that
\begin{displaymath}
\mu(\tilde{\mathcal{D}}_i)=1-\langle c_1(\spinc_z(\{x_1,\bullet\})),[F]\rangle .
\end{displaymath}
In fact
\begin{equation}
\spinc_z(\{x_1\}\cup \mathfrak{x})=\spinc(\{x\}\cup\mathfrak{x}),
\end{equation}
where the right hand side is the $\SpinC$ structure over $S^3_0(K)$ assigned
to the Kauffman state $\{x\}\cup\mathfrak{x}$ in the Ozsv\'ath-Szab\'o
Floer theory.

Suppose that the Ozsv\'ath-Szab\'o Floer homology $\widehat{\HFK}(K)$ is non-zero in
the $\SpinC$ structure $s$, and that $s$ is the highest $\SpinC$ structure
with this property.
By the result of Ozsv\'ath and Szab\'o in \cite{OS-genus}, this is the same as assuming that
$s$ is the genus of the knot $K$.
Let $\{x\}\cup\mathfrak{x}_1,...,\{x\}\cup\mathfrak{x}_k$
be the Kauffman states of $K$ in the $\SpinC$ class $s$; and that
$\{x\}\cup\mathfrak{y}_1,...,\{x\}\cup\mathfrak{y}_l$ are the Kauffman states
in the higher $\SpinC$ classes, say $\spinc(\{x\}\cup\mathfrak{y}_j)=s+i_j-1$.
Here the integer $s$ represents a $\SpinC$ structure in $\SpinC(S^3_0(K))$ via
the isomorphism $\SpinC(S^3_0(K))\simeq \Z$. 

The disk between $\{x\}\cup\mathfrak{y}_j$ and $\{x\}\cup\mathfrak{y}_i$, if they
are in the same $\SpinC$ class, is the same as the disk between
$\{x_r\}\cup\mathfrak{y}_j$ and $\{x_r\}\cup\mathfrak{y}_i$.

Look at the $\SpinC$ structure $t$ assigned by the map $\spinc_z$ to
$\{x_1\}\cup\mathfrak{x}_i$ and let $\ut=t-\tfrac{1}{2}$. The Kauffman
states in this $\SpinC$ structure are exactly the following:
\begin{displaymath}
\begin{CD}
\{x_1\}\cup\mathfrak{x}_i,\{y_1\}\cup\mathfrak{x}_i,\ \ \ \ \ \ i=1,...,k,\\
\{x_{i_j}\}\cup\mathfrak{y}_j,\{y_{i_j}\}\cup\mathfrak{y}_j,\ \ \ \ \ \ j=1,...,l.\\
\end{CD}
\end{displaymath}
First consider the disks supported on the part of the surface outside the cylinder.
The Kauffman states of the second type cancel each other since the disks
between them are in fact identical to the disks in the higher $\SpinC$ structures
of $\widehat{\HFK}(K)$, which give trivial groups. This may be thought of as puncturing
the domains on the cylinder and doing the cancellations via computing the homology
of the resulting Heegaard diagram. 

Among the Kauffman states of the form $\{x_1\}\cup\mathfrak{x}_i$ the cancellations
are also identical to the cancellations of the Ozsv\'ath-Szab\'o hat theory.
In fact, the main possible problem is the possibility of a set of boundary maps
of the form:
\begin{displaymath}
\begin{array}{ccc}
\{x_{i_j}\}\cup\mathfrak{y}_j& & \{x_1\}\cup\mathfrak{x}_i\\
\downarrow & \stackrel{\pi}{\searrow} &\downarrow\\
\{x_1\}\cup\mathfrak{x}_{i'} & &\{x_{i_m}\}\cup\mathfrak{y}_m,
\end{array}
\end{displaymath}
where $\pi$ is coming from a disk supported away from the cylinder as
above, imposing a cancellation of $\{x_{i_j}\}\cup\mathfrak{y}_j$
against $\{x_{i_m}\}\cup\mathfrak{y}_m$.
This cannot happen, since the map to $\{x_1\}\cup\mathfrak{x}_i$ or
the map from $\{x_1\}\cup\mathfrak{x}_i$ has very negative coefficients
in some domains in the cylinder.

Similarly, the cancellations among $\{y_1\}\cup\mathfrak{x}_i$'s are identical
to those in the hat theory.

To understand $\widehat{\HFL}(K,\ut)$, we should study the boundary maps between
$\{x_1\}\cup\mathfrak{x}_i$ and $\{y_1\}\cup\mathfrak{x}_j$. Note that
the domain of any disk from $\{x_1\}\cup\mathfrak{x}_i$ to $\{y_1\}\cup\mathfrak{x}_j$
has negative coefficients in the cylinder. Thus there is no boundary map in this direction.

Potentially there can be a  boundary map from $\{y_1\}\cup\mathfrak{x}_j$ to
$\{x_1\}\cup\mathfrak{x}_i$. Let $\mathcal{D}'$ denote the domain
of the disk from $\{x\}\cup\mathfrak{x}_j$ to $\{x\}\cup\mathfrak{x}_i$ which is supported
outside the cylinder. Then the domain of the disk from $\{y_1\}\cup\mathfrak{x}_j$ to
$\{x_1\}\cup\mathfrak{x}_i$ will  be equal to
$\mathcal{D}'+\tilde{\mathcal{D}}_1$, and the Maslov index is
\begin{equation}
\begin{split}
\mu(\mathcal{D}')+\mu(\tilde{\mathcal{D}}_1)&=\mu(\mathfrak{x}_j)-\mu(\mathfrak{x}_i)
+\mu(\tilde{\mathcal{D}}_1)\\
&=\mu(\mathfrak{x}_j)-\mu(\mathfrak{x}_i)+1-\langle c_1(s),[F]\rangle .
\end{split}
\end{equation}
Since $s$ is  the highest nontrivial $\SpinC$ structure for which
$\widehat{\HFK}(K,s)$ is nonzero, $\langle c_1(s),[F]\rangle $ is
at least $2$ (unless $K$ is the trivial knot). To have a disk from
$\{y_1\}\cup\mathfrak{x}_j$ to $\{x_1\}\cup\mathfrak{x}_i$ we need
to have $\mu(\mathfrak{x}_j)-\mu(\mathfrak{x}_i)\geq 2$. If
$\mathfrak{x}_1$ is the Kauffman state with highest Maslov grading
among $\mathfrak{x}_j$s which survives the cancellations in the
standard hat theory, then this condition may not be satisfied.
Thus $\{y_1\}\cup\mathfrak{x}_1$ will not be cancelled at the
level of homology. The result is the following
theorem:

\begin{thm}
Suppose that $K$ is a nontrivial  knot in $S^3$. Then the longitude Floer homology $\widehat{\HFL}(K)$
is nontrivial.\end{thm}

The $\SpinC$ structure determined by $\spinc_z(\{x_1\}\cup \mathfrak{x}_j)$ may be described
as $s=\spinc(\{x\}\cup \mathfrak{x}_j)$. For the $\SpinC$ structures $t>s$ the above
argument shows that in fact the longitude Floer homology is trivial.
Thus the element of
$$\tfrac{1}{2}\text{PD}[\mu]+\SpinC(S^3_0(K))\simeq \tfrac{1}{2}+\Z$$
associated with $\{x_1\}\cup \mathfrak{x}_j$ via the map
$\relspinc$ is $\spinc(\{x\}\cup\mathfrak{x}_j)-\frac{1}{2}$.

We may do the winding in the other direction. This time a similar argument shows
that there is a minimum $\SpinC$ structure described as $\us'=s(\{x\}\cup\mathfrak{x}_j')+\frac{1}{2}$
such that for $\ut<\us'$ the longitude Floer homology is trivial and it is nontrivial
for $\us'$. Here $\mathfrak{x}_j'$s are the Kauffman states in the usual Heegaard
diagram of $K$ which produce the lowest nontrivial group $\widehat{\HFK}(K,s'-\frac{1}{2})$.

Possibility of winding in the two different directions and the symmetry of
$\widehat{\HFK}(K)$ implies the
existence of a symmetry in the
longitude Floer homology. In fact we may prove the following theorem:

\begin{thm}\label{thm:genusbound} Suppose that $g$ is the genus of a nontrivial knot $K$ is $S^3$. Then
$$
\widehat{\HFL}(K,g-\tfrac{1}{2})\simeq \widehat{\HFL}(K,-g+\tfrac{1}{2})\neq 0,$$
and for any $\ut>g$ in $\tfrac{1}{2}+\mathbb{Z}$ the groups
$\widehat{\HFL}(K,\ut)$ and $\widehat{\HFL}(K,-\ut)$ are both trivial. Furthermore,
for any $\ut$ in $\tfrac{1}{2}+\mathbb{Z}$ there is an isomorphism of the
relatively graded groups (graded by the Maslov grading),
$$
\widehat{\HFL}(K,\ut)\simeq \widehat{\HFL}(K,-\ut).
$$\end{thm}

By a result of Ozsv\'ath and Szab\'o (\cite{OS-fibered}), we know that for a fibered knot $K$ of genus $g>0$,
there exists a Heegaard diagram with a single generator in highest $\SpinC$ structure which is $s=g$.
Furthermore for the $\SpinC$ structures $s>g$, there is no other generator of this Heegaard diagram.
Using this Heegaard diagram in the above argument we obtain a Heegaard diagram for the longitude Floer homology
with two generators in the $\SpinC$ structure $\us=g-\frac{1}{2}$ and no generators in the $\SpinC$ structures
$\us>g$. Furthermore, the above argument shows that the two generators
in this $\SpinC$-structure can not cancel each other because of the
difference in their Maslov gradings. As a result we obtain the following:

\begin{prop} If $K$ is a fibered knot of genus $g$ then
$\widehat{\HFL}(K,\pm(g-\frac{1}{2}))$ is equal to
$\mathbb{Z}\oplus{\mathbb{Z}}$.\end{prop}

\section{Example: $T(2,2n+1)$}

We continue by an explicit computation of the longitude Floer
homology for the $(2,2n+1)$ torus knots.
\begin{figure}[ht!]\small\anchor{fig:longitude.3}
\psfrag {DL}{$D_L$}
\psfrag {DR}{$D_R$}
\psfrag {i}{$i^{\rm th}$ intersection}
\psfrag {p}{point}
\cl{\hspace{-2cm}\includegraphics{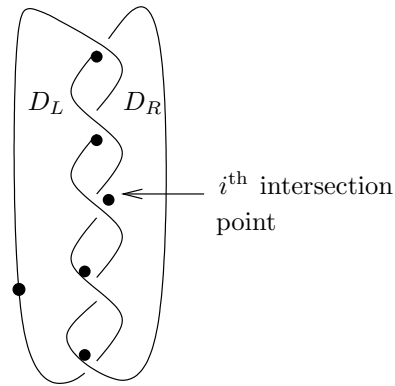}}
\caption{\label{fig:longitude.3}
{The Kauffman state $\mathfrak{z}_i$ of the torus knot is shown.  Here
the torus knot is the $(2,5)$-knot and the Kauffman state is $\mathfrak{z}_3$. }}
\end{figure}

We remind the reader that associated with any planar diagram for a
knot $K$, and a marked point on it, is a Heegaard diagram for the knot
$K$ as discussed in \cite{OS-alternating}. The generators of the complex
(i.e.\ Heegaard Floer complex defined by Ozsv\'ath-Szab\'o \cite{OS-knot}
and Rasmussen \cite{Ras2})
associated with this Heegaard diagram may be described as follows.
If $A_1,A_2,\ldots,A_m$ are the regions in the complement of the the knot
in plane which are not neighbors of the marked point, then any
generator corresponds with an $m$-tuple of points such that in each
region $A_i$ exactly one marked point is chosen. Each marked point is
located near a self intersection in the plane projection of the knot
(obtained from the knot diagram). Furthermore, for any such $m$-tuple,
it is required that from the four quadrants in each self intersection,
exactly in one of them a marked point is chosen. These sets of marked
points are called the \emph{Kauffman states} for the planar knot
diagram. If the unbounded region is a neighbor of the marked point on
the diagram (which is the case in what follows), the meridian curve
will intersect a unique $\alpha$-curve in a single point, and any
generator will contain this intersection point.

Consider a standard plane diagram of the $(2,2n+1)$ torus knot
shown in   \figref{fig:longitude.3}.  Let the bold points denote the
marked points representing a Kauffman state. Other than the outside region,
there are two large bounded regions on the right side and the left side of
the twists, denoted by $D_R$ and $D_L$ respectively, and $2n$ small regions,
which we call $D_1,\ldots,D_{2n}$ from top to bottom. There is a marked point
on the knot which is put in the common boundary of $D_L$ and the unbounded region.

\begin{figure}[ht!]\small\anchor{fig:example}
\psfrag {l}{$l$}
\psfrag {a}{$\alpha_1$}
\psfrag {x1}{$x_1$}
\psfrag {y1}{$y_1$}
\psfrag {x2}{$x_2$}
\psfrag {y2}{$y_2$}
\cl{\hspace{-2cm}\includegraphics{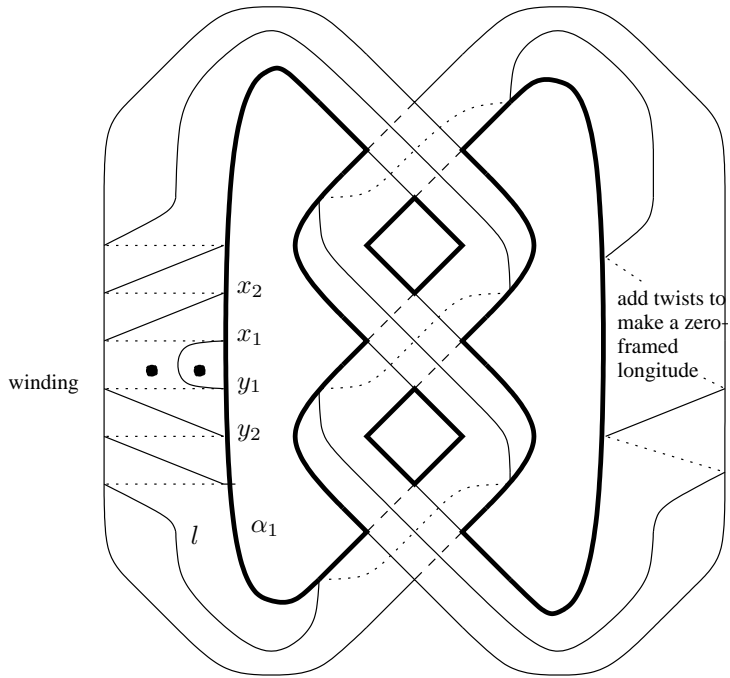}}
\caption{\label{fig:example}
{The Heegaard diagram associated with the trefoil is presented.
On the handle appearing on the right-hand-side we may do enough twists
so that the diagram represents a three-manifold with $b_1=1$.  The
winding is done on the left-hand-side handle.  The bold curves are the
$\alpha$ curves and the rest of them are the $\beta$ curves.
}}
\end{figure}

Let $\mathfrak{z}_i$ be the Kauffman state containing a marked point in $D_R$ at the $i$-th
intersection. There is a unique Kauffman state described by this property.
Moreover, the states $\mathfrak{z}_1,\ldots,\mathfrak{z}_{2n+1}$ will be all of the
Kauffman states of (generators of the Heegaard Floer complex for) the
$(2,2n+1)$ torus knot. As it was noted earlier the Kauffman states are
in one-to-one correspondence with the generators. So each
$\mathfrak{z}_i$ may be thought of  as a set of $2n+1$ intersection
points in the Heegaard diagram which, together with the unique
intersection point on the meridian, give a generator for the Heegaard
Floer complex. These two alternative ways of thinking about the
Kauffman states $\mathfrak{z}_i$ are used in the following.

The $\SpinC$ grading of the Kauffman states is described via
$\spinc(\mathfrak{z}_i)=i-n-1$,
and the (relative) Maslov grading by $\mu(\mathfrak{z}_i)=i-1$, all in
the sense of \cite{OS-knot}. Note that $\spinc(\mathfrak{z}_i)\in \Z=
\SpinC(S^3_0(K))$ is the well-defined $\SpinC$ structure used
in the Heegaard Floer homology of Ozsv\'ath-Szab\'o (\cite{OS-knot})
and Rasmussen (\cite{Ras2}).

 After winding $l$ along $m$ sufficiently many
 times, the proof of theorem~\ref{thm:genusbound} (cf.\ section 5 of \cite{OS-3m1})
 may be copied to prove the following:

 \begin{lem}
 For any $\SpinC$ class
$$\us\in \tfrac{1}{2}+\Z=\tfrac{1}{2}\text{PD}[\mu]+\SpinC(S^3_0(K)),$$
 with the property $|\us|<n=\text{genus}(T(2,2n+1))$,
 all the generators in the class $\us$ are of the form:
$$
 \mathfrak{x}_{ij}=\{x_i\}\cup \mathfrak{z}_j,\ \ \ \ \
 \mathfrak{y}_{ij}=\{y_i\}\cup \mathfrak{z}_j,
 $$
where $\mathfrak{z}_j$s are considered as sets of $2n+1$ intersection
points in the Heegaard diagram, and $x_1,x_2,\ldots,y_1,y_2,\ldots$ are the
intersection points on $l$ which result from winding it around the
meridian $\mu$.
 \end{lem}
These generators will be called
\emph{Kauffman states for longitude Floer homology} or just \emph{Kauffman
states} if it is clear from the context that longitude Floer complex is
considered.

 It is easy to check that the following assignments, satisfy all the relative
 Maslov grading computations and the equations for $\SpinC$ differences:
 \begin{equation}
 \begin{CD}
 -\mu(\mathfrak{y}_{ij})=\mu(\mathfrak{x}_{ij})=j-n-\frac{3}{2},\\
  \relspinc(\mathfrak{x}_{ij})=\relspinc(\mathfrak{y}_{ij})=j-i-n-\frac{1}{2}.
 \end{CD}
 \end{equation}
Here we are assigning rational values as the Maslov grading, which is
an abuse of notation. However, note that here we are only interested
in relative grading, and the relative Maslov gradings are still by integers.

The Kauffman states which lie in the $\SpinC$ structure
$\us=s-\frac{1}{2}\in \frac{1}{2}+\mathbb{Z}$
are those $\mathfrak{x}_{ij}$ and $\mathfrak{y}_{ij}$ for which
$j-i=n+s$.

Remember that there cannot be any boundary map going from $\mathfrak{x}_{ij}$
to $\mathfrak{y}_{ij}$. Furthermore, if there
exists a map from $\mathfrak{x}_{ij}$ to $\mathfrak{x}_{i+k,j-k}$ regardless of what
$N$ is, then there is a map from $\mathfrak{x}_{l,j}$ to $\mathfrak{x}_{l+k,j-k}$
regardless of what $N$ is, for all other $l$. Conversely, if there is no map from $\mathfrak{x}_{ij}$ to $\mathfrak{x}_{i+k,j-k}$ regardless of what
$N$ is, then there is no map from $\mathfrak{x}_{l,j}$ to $\mathfrak{x}_{l+k,j-k}$.
This is because of the isomorphism between the domains of the disks between the corresponding
generators.

Since we already know that $\CFLhat(K)$ is symmetric with respect to the $\SpinC$ structure
$$\us\in \tfrac{1}{2}+\Z\simeq \tfrac{1}{2}\text{PD}[\mu]+\SpinC(S^3_0(K)),$$
and since the genus of $K=T(2,2n+1)$ is $n$,
it is enough to compute $\CFLhat(K,\us)$ for the $\SpinC$ structures $0<\us<n$.
Here $\mu$ represents the meridian of the knot $K$.

If $0<\us<n$, then because of the above Maslov grading of the
 generators, there
cannot be any boundary maps from any of $\mathfrak{y}_{ij}$'s
to any of $\mathfrak{x}_{kl}$'s in $\SpinC$
class $\us$. Thus, the only boundary maps that should be studied
are the boundary maps within $\mathfrak{x}_{ij}$'s,
as well as the boundary maps
within $\mathfrak{y}_{kl}$'s.

\begin{figure}[ht!]\small\anchor{fig:Tdomain}
\psfrag {a}{(a)}
\psfrag {b}{(b)}
\cl{\hspace{-2cm}\includegraphics[width=3.5in]{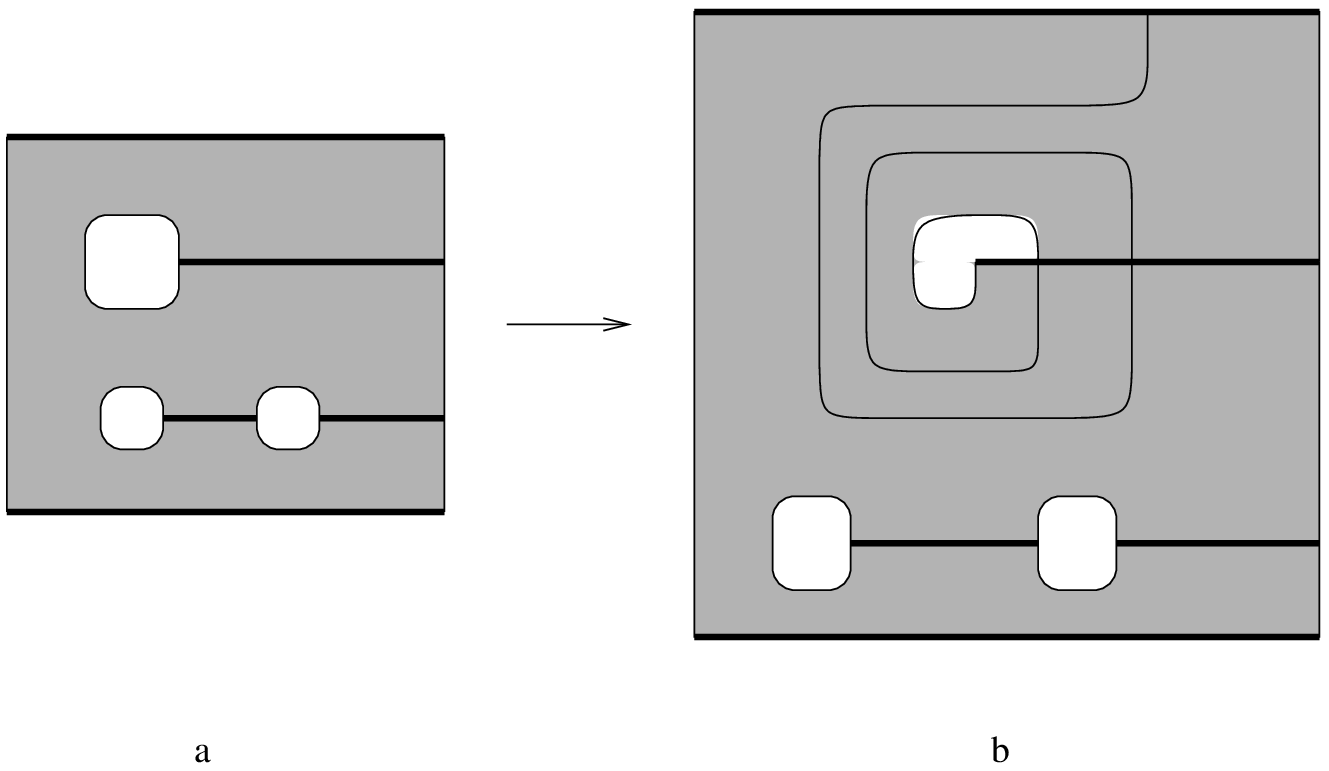}}
\caption{\label{fig:Tdomain}
{The domain between $\mathfrak{x}_{ij}$ and
$\mathfrak{x}_{(i-1)(j-1)}$ (b)
is a modification of the domain
connecting the two generators $\mathfrak{z}_j$ and
$\mathfrak{z}_{j-1}$ (a)
in the
Heegaard diagram obtained from the alternating projection of $K$.
If there are $k$ small circles in the domain on the left, we will
denote it by $\sD_k$.   Note that the bold curves are in $\alphas$
while the regular curves are in $\betas$.
}}
\end{figure}

If $\mathfrak{x}_{kl}$ appears in the boundary of $\mathfrak{x}_{ij}$, then they are
in the same $\SpinC$ class and the Maslov grading of the first generator is
one less than the Maslov grading
of the second generator. This implies that $k=i-1$ and $j=l-1$. The domain between
$\mathfrak{x}_{ij}$ and $\mathfrak{x}_{(i-1)(j-1)}$ is a modification of the domain
connecting the two generators $\mathfrak{z}_j$ and $\mathfrak{z}_{j-1}$ in the
Heegaard diagram obtained from the alternating projection of $K$ shown in
\figref{fig:longitude.3}. If $j=2l+1$ for some $l$, then the domain between
$\mathfrak{z}_j$ and $\mathfrak{z}_{j-1}$ is illustrated in \figref{fig:Tdomain}(a),
while the modified domain connecting $\mathfrak{x}_{ij}$ and $\mathfrak{x}_{(i-1)(j-1)}$
will be of the type shown in \figref{fig:Tdomain}(b). Let us denote this
last domain by
$\sD_k$, where $k$ is the number of circles inside the rectangle. 

The moduli spaces $\sM(\sD_k)$ and $\sM(\sD_{k-1})$ are in fact
cobordant, since $\sD_k$ is obtained from $\sD_{k-1}$ via the
operation of adding a handle. This can be proved using the usual
argument of Ozsv\'ath and Szab\'o for the invariance of the Floer
homology when we add a one handle to the surface, and a pair of
cancelling
curves to $\alphas$ and $\betas$ (see \cite{OS-3m1}).

To show that the total contribution of the domain $\sD_k$ to the boundary map is
$\pm 1$, we only have to show this for $\sD=\sD_0$.

\begin{lem}
Let $\sD$ be as above. Then the algebraic sum of the points in the moduli space
$$\widehat{\sM}(\sD)=\frac{\sM(\sD)}{\R}$$
is $\pm 1$.
\end{lem}

\begin{proof}
Consider the embedding of the domain $\sD$ in a genus three Heegaard diagram
which is shown in \figref{fig:Tdomain2}. Let the bold and the regular curves
denote the $\alpha$ and the $\beta$ curves respectively. Choose $\alpha_i$'s and
$\beta_j$'s so that the $\alpha$ curve which spins around the center is $\alpha_1$
and the $\beta$ curve cutting it several times is $\beta_1$.

Consider the small dotted circle $\theta_1$ in
\figref{fig:Tdomain2} and complete it into a set of three
disjoint linearly independent simple closed curves by adding
Hamiltonian isotopes of the curves $\beta_2$ and $\beta_3$, which
we will call $\theta_2$ and $\theta_3$ respectively. We choose
them so that $\theta_i$ intersects $\beta_i$ in a pair of
transverse cancelling intersection points, for $i=2,3$. Call the
resulting
sets of curves $\alphas, \betas$ and $\thetas$.

The triple Heegaard diagram
$$H=\Big\{\Sig_3,\alphas,\thetas,\betas;u,v,w\Big\},$$
with $u,v$ and $w$ being the marked points of \figref{fig:Tdomain2}, induces a chain
map
$$\sF:\widehat{CF}(\alphas,\thetas)\otimes \widehat{CF}(\thetas,\betas)\lra
\widehat{CF}(\alphas,\betas).$$
The map $\sF$ is defined through a count of holomorphic triangles which miss
the marked points $u,v$ and $w$ (see \cite{OS-3m1,OS-3m2} for more details on the
construction of $\sF$). The complex $\widehat{CF}(\thetas,\betas)$
gives the Floer homology
associated with the three-manifold $(S^1\times S^2)\#(S^1\times S^2)$. There is a
top generator of this homology group which we may denote by $\Theta$. The complex
$\widehat{CF}(\alphas,\thetas)$
has precisely two generators $\x$ and $\y$, with a single
boundary map going from $\x$ to $\y$. The image
$\sF(\x \times \Theta)$ will have several terms, probably in different $\SpinC$
classes.

\begin{figure}[ht!]\small\anchor{fig:Tdomain2}
\psfraga <0pt, 2pt> {u}{$u$}
\psfraga <0pt, 2pt> {v}{$v$}
\psfraga <0pt, 2pt> {w}{$w$}
\psfrag {t}{$\theta_1$}
\psfrag {a1}{$\alpha_1$}
\psfrag {a2}{$\alpha_2$}
\psfrag {a3}{$\alpha_3$}
\psfrag {b1}{$\beta_1$}
\psfrag {b2}{$\beta_2$}
\psfrag {b3}{$\beta_3$}
\cl{\hspace{-1.5cm}\includegraphics{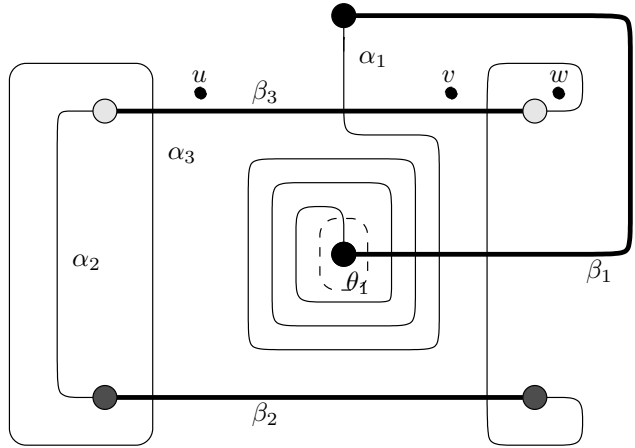}} \caption{\label{fig:Tdomain2}
{The domain $\sD$ may be embedded in a genus three Heegaard
diagram. The curve winding around the center is $\alpha_1$, which
is completed to $\alphas=\{\alpha_1,\alpha_2,\alpha_3\}$. The
curve $\beta_1\in \betas=\{\beta_1,\beta_2,\beta_3\}$ cuts
$\alpha_1$ several times. The dotted small circle is $\theta_1$
which is completed to a triple $\thetas$ by adding the Hamiltonian
isotopes $\theta_2$ and $\theta_3$ of $\beta_2$ and $\beta_3$. The
intersection points between $\beta_1$ and $\alpha_1$ are labelled
$x_1,x_2,\ldots$ with $x_1$ the intersection point on the right hand
side of $\theta_1$ in the picture. }}
\end{figure}

Denote the intersection points of $\alpha_1$ and $\beta_1$ in the spiral by
$x_1,x_2,\ldots$, so that $x_1$ is the one that is closest to the center of the spiral.
Each $x_i$ may be completed to a generator of
$\widehat{CF}(\alphas,\betas)$ precisely in two ways,
which will be denoted  by $\{x_i\}\cup \z$ and $\{x_i\}\cup \w$. We may choose them so
that $\{x_1\}\cup \z$ and $\{x_2\}\cup\w$ are in the same $\SpinC$ class. Under this
assumption the domain connecting them is the domain
$\sD$ introduced above. Denote this same $\SpinC$
class by $s$, and denote the part of the image of $\sF$ in the $\SpinC$ class $s$ by
$\sF_s$. Clearly $\sF_s$ is also a chain map.

It is not hard to check, using the energy filtration of \cite{OS-3m2}, that
we would have
$$\sF_s(\x\otimes \Theta)=\pm \{x_1\}\cup \z + \text{lower energy terms, and}$$
$$\sF_s(\y\otimes \Theta)=\pm \{x_2\}\cup \w +\text{lower energy terms}.$$
It is then an algebraic fact that $\{x_2\}\cup \w$ appears in the boundary of
$\{x_1\}\cup \z$ with coefficient $\pm 1$, which is on its own a result of
$\partial(\x \otimes \Theta)=\y \otimes \Theta$.
This completes the proof of the lemma.
\end{proof}

This lemma implies that $\partial(\mathfrak{x}_{(i)(2l+1)})=
\pm \mathfrak{x}_{(i-1)(2l)}$.
Since $\partial \circ \partial=0$, we may conclude that
$\partial (\mathfrak{x}_{(i)(2l)})=0$ for all $i,l$, unless $i$ is too large
(i.e.\ irrelevant). Similarly we may deduce that
$$\partial(\mathfrak{y}_{(i)(2l)})=\mathfrak{y}_{(i-1)(2l-1)},\text{ and}$$
$$\partial(\mathfrak{y}_{(i)(2l+1)})=0.$$
We may summarize all these information as the following theorem.

\begin{thm}
Let $K=T(2,2n+1)$ denotes the $(2,2n+1)$ torus knot in $S^3$. Then for any $\SpinC$
structure $$\us\in\tfrac{1}{2}+\Z\simeq \tfrac{1}{2}\text{PD}[\mu]+\SpinC(S^3_0(K)),$$
the longitude Floer homology of $K$ is trivial if $|\us|>n$. Otherwise it is given by
$$\HFLhat(K,\us)=\Z_{(-n+\frac{1}{2})}\oplus \Z_{(\epsilon(\us)\us)},$$
where $\epsilon(\us)=(-1)^{n-\frac{1}{2}-\us}$, and $\Z_{(k)}$ denotes a copy of $\Z$ in
(relative) Maslov grading $k$.
\end{thm}
\begin{proof}
The proof for $\us>0$ is just an algebraic result of the cancellations induced by the
map $\partial$ above. For $\us<0$, it is the result of the symmetry on $\HFLhat(K)$.
\end{proof}

\section{A Heegaard diagram for Whitehead double}
In this section, we will construct an appropriate Heegaard diagram for $K_L$ out of a Heegaard
diagram for $K$.

\begin{figure}[p]\small\let\tiny\footnotesize\anchor{fig:whitehead.2}
\psfrag {a}{(a)}
\psfraga <20pt,0pt> {b}{(b)}
\psfrag {m}{$m$}
\psfrag {l}{$l$}
\psfrag {d}{$\delta_1$}
\psfraga <0pt,-2pt> {n}{\tiny$n$}
\psfrag {g}{\tiny$\gamma$}
\psfraga <2pt,0pt> {z}{\tiny$z$}
\psfrag {w}{\tiny$w$}
\psfrag {la}{\tiny$\lambda$}
\psfrag {a1}{\tiny$\alpha_1$}
\psfrag {a2}{\tiny$\alpha_2$}
\psfrag {a3}{\tiny$\alpha_3$}
\psfrag {a4}{\tiny$\alpha_4$}
\psfraga <0pt,-2pt> {b1}{\tiny$\beta_1$}
\psfraga <0pt,-2pt> {b2}{\tiny$\beta_2$}
\psfraga <0pt,-2pt> {b3}{\tiny$\beta_3$}
\psfraga <0pt,-2pt> {b4}{\tiny$\beta_4$}
\cl{\hspace{-2cm}\includegraphics[width=4.4in]{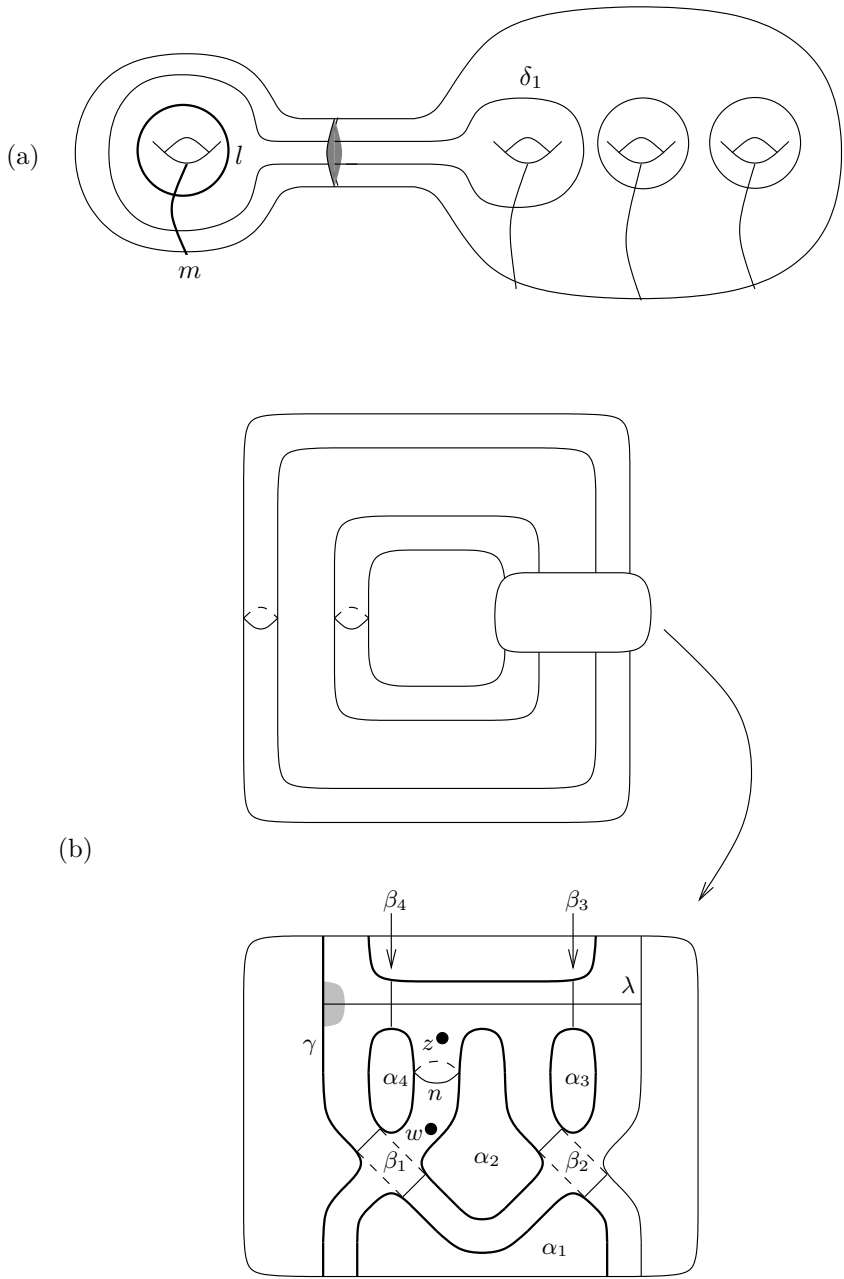}}
\caption{\label{fig:whitehead.2} {(a) A Heegaard diagram
associated with $K$. Here $m$ denotes the meridian, and $l$ is the
longitude of the knot. The curve $\delta_1$ is the unique $\delta$
curve cutting $m$. (b) A Heegaard diagram for $L$ in the solid
torus. The thicker curves denote the $\alpha$ curves, and the
thinner ones are $\beta$'s. The curve $n$ denotes the meridian of
$L$. }}
\end{figure}

 Suppose that
$(\Sigma;\deltas;\{m=\gamma_1\}\cup \gammas_0)$
is a Heegaard diagram for $K$ together with an extra curve $l$ (As usual,
$\deltas=\{\delta_1,\ldots,\delta_g\}$ and $\gammas_0=\{\gamma_2,\ldots,\gamma_g\}$).
Here $l$ is the curve
with the property that it
intersects the meridian $m=\gamma_1$ exactly once but does not cut any
other $\gamma$ curve. The curve
$l$ represents the longitude of the knot $K$ in such a way that
$$(\Sigma;\deltas;\{ l\}\cup \gammas_0)
$$is a Heegaard diagram for $S^3_0(K)$.
Define $$\gammas=\{l\}\cup \gammas_0.$$
 We may bring this Heegaard diagram into the form
shown in \figref{fig:whitehead.2}, where the surface $\Sigma$ is a connected sum
$\Sigma= T\# S$ of the  torus $T=S^1\times S^1$ and $S$; a surface of genus
$g-1$. We assume that $m$ and $l$ are the standard generators of the homology
of $T$, and that all other $\gamma$ curves are on $S$ and are in the standard
configuration such that by attaching a disk to these $\gamma$ curves we get the
handlebody formed by the inside of $S$. The diagram
$(\Sigma,\deltas;\gammas_0)$
gives a Heegaard diagram for the complement of the knot $K$ in $S^3$ and the
neighborhood of $K$ may be identified with the interior of $T$, since
$l$ has zero linking number with the core of the torus $T=S^1\times S^1$.

If we embed the knot $L$ inside the torus $T$ and find a Heegaard diagram for
its complement, this Heegaard diagram together with the Heegaard diagram
for $K$ will give a diagram for the double of the knot in $S^3$.

More precisely, consider the Heegaard diagram shown in \figref{fig:whitehead.2}~(b) for
the unknot sitting inside the solid torus. Here the thick curves denote the $\alpha$
circles, while the thin ones are $\beta$'s. There is an extra curve $\lambda$
shown in the picture, which we save for the later purposes. There is a special
$\alpha$-curve denoted by $\gamma$ in the figure, which represents the generator
of the first homology group $H_1(D^2\times S^1,\mathbb{Z})$ of the solid torus.

If we attach a disk to each $\alpha$ curve, except for $\gamma$,
in this Heegaard diagram and a disk to each
of the $\beta$ circles other than the meridian $n$, we will get the complement of
the knot $L$ inside a solid torus $S^1\times D^2$.

Put this solid torus inside the torus $T$. Attach the surface of the solid
torus and $T$ by a one-handle connecting the intersection of $\gamma_1$ and $l$ on $T$ to
the intersection of $\gamma$ and $\lambda$ on the Heegaard diagram for $L$.

The result of this operation may be regarded as a connected sum of the surfaces
$\Sigma$ and $C$. Here $(\Sigma,\deltas,\{m\}\cup\gammas_0)$ is the above Heegaard diagram for $K$, and
$C$ is the surface in the Heegaard diagram of $L$ used above.

Denote by $(n,\beta_1,\ldots,\beta_4)$ the $\beta$ curves on $C$ and by
$(\gamma,\alpha_1,\ldots,\alpha_4)$ the $\alpha$ curves, as is shown in \figref{fig:whitehead.2}~(b).

\begin{figure}[ht!]\small\anchor{fig:whitehead.3}
\psfrag {m}{$m$}
\psfrag {l}{$l$}
\psfrag {a}{$\alpha$}
\psfrag {la}{$\lambda$}
\psfrag {g}{$\gamma$}
\psfraga <-37pt,-10pt> {cuta}{\vbox{cut along the plane\newline
passing though $l,\alpha$}}
\psfraga <5pt,-10pt>  {cutm}{cut along the curve $m$}
\cl{\hspace{-1.5cm}\includegraphics[width=3.9in]{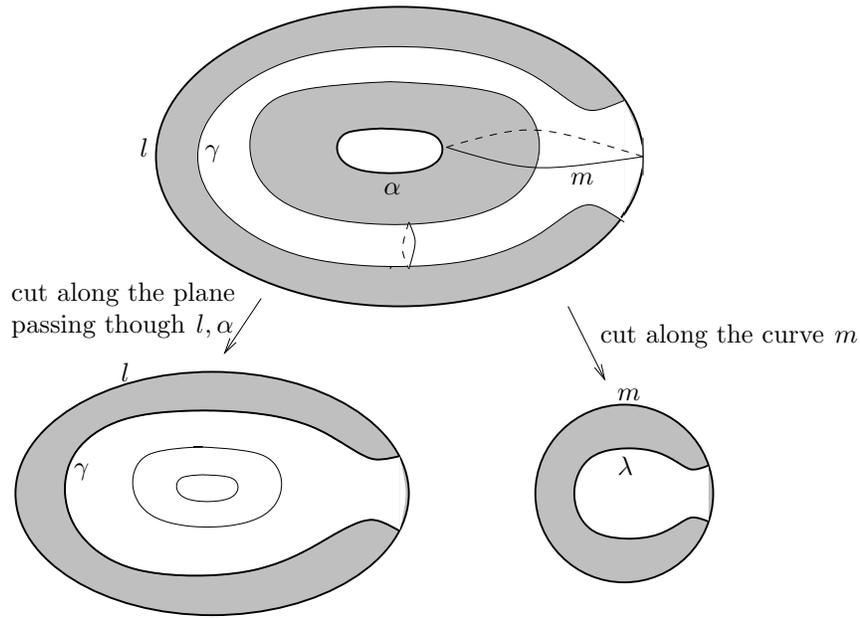}}
\caption{\label{fig:whitehead.3}
{If we cut the torus $T$ by a horizontal plane, the intersection will be as
shown on the left side. There is a disk which is dashed in this picture
with boundary $\gamma\# l$. If the cut is vertical and on the connecting handle,
the picture is as shown on the right. Again, there is a disk with boundary
$\lambda\# m$}}.
\end{figure}

In order to find a Heegaard diagram for the complement of the Whitehead double
in the sphere $S^3$, we have to fill the space between the solid torus and the
torus $T$. Looking at $T\#C$, there are two disks which sit in the empty space
between $T$ and the solid torus. Namely, if we cut the union by a horizontal plane,
the intersection will look like the left hand side of \figref{fig:whitehead.3}.
There is a disk bounded by the connected sum $\gamma\# l$. This disk is
dashed in \figref{fig:whitehead.3}.

We may also cut the torus with a vertical plane. If the cut is made in a way that it passes
through $m$ and $\lambda$,
and cuts the handle connecting the
solid torus and the torus $T$,  then the cut will look like what is shown on the right
hand side of \figref{fig:whitehead.3}. Again, there is a disk which is dashed in
the picture,
with a  boundary which is the connected sum $\lambda \#m=\lambda\#\gamma_1$.

The result of this operation is a Heegaard diagram for the Whitehead double of $K$:
$$
(\Sigma\#C;\{n,\beta_1,\ldots,\beta_4\}\cup \deltas;\{\alpha_1,\ldots,\alpha_4\}\cup
\{\lambda\#m,\gamma\# l\}\cup \gammas_0;z,w),
$$
where $z$ and $w$ are two base points which are put on the two sides of the curve $n$ on $C$.
We will use this Heegaard diagram to relate the Ozsv\'ath-Szab\'o Floer homology
of the Whitehead double of $K$ to the longitude Floer homology discussed in the
earlier sections.

\section{Whitehead double; homology computation}
In order to obtain the Ozsv\'ath-Szab\'o Floer homology groups we should first form the
chain complex by identifying the relevant generators of this Heegaard diagram.

There are two types of generators for this Heegaard diagram:

(1)\qua The Kauffman states which are in correspondence with a pair of
generators of the form $\{\x,\y\}$, where $\x$
is a generator of $$(C;\gamma,\alpha_1,\ldots,\alpha_4;n,\beta_1,\ldots,\beta_4),$$ and
$\y$ is a generator of $(\Sigma,\deltas;\{m\}\cup \gammas_0)$.
We call these generators \emph{meridian Kauffman states}.

(2)\qua The Kauffman states which are associated with a pair of generators of the form $\{\x,\y\}$, where $\x$
is a generator of $$(C;\lambda,\alpha_1,\ldots,\alpha_4;n,\beta_1,\ldots,\beta_4),$$ and
$\y$ is a generator of $(\Sigma,\deltas;\{l\}\cup \gammas_0)$.
We call these generators the \emph{longitude Kauffman states}.

\begin{figure}[ht!]\small\anchor{fig:whitehead.4}
\psfraga <-2pt,0pt>  {-M}{$-M$}
\psfraga <-2pt,0pt>  {M}{$M$}
\psfraga <-2pt,0pt>  {1}{$(1)$}
\psfraga <-2pt,0pt>  {0}{$0$}
\psfraga <-2pt,0pt> {a}{\tiny$a$}
\psfraga <-2pt,0pt>  {b}{\tiny$b$}
\psfraga <-2pt,0pt>  {c}{\tiny$c$}
\psfraga <-2pt,2pt>  {d}{\tiny$d$}
\psfraga <-2pt,1pt>  {e}{\tiny$e$}
\psfraga <-2pt,0pt>  {f}{\tiny$f$}
\psfraga <-2pt,0pt>  {g}{\tiny$g$}
\psfraga <-2pt,0pt>  {h}{\tiny$h$}
\psfraga <-2pt,0pt>  {k}{\tiny$k$}
\psfraga <-2pt,0pt>  {l}{\tiny$l$}
\psfraga <-2pt,0pt>  {m}{\tiny$m$}
\psfraga <-2pt,0pt>  {n}{\tiny$n$}
\psfraga <-2pt,0pt>  {p}{\tiny$p$}
\psfraga <-2pt,1pt>  {q}{\tiny$q$}
\psfrag {r}{\tiny$r$}
\psfraga <1pt,2pt>  {a'}{\tiny$a'$}
\psfraga <2pt,2pt>  {b'}{\tiny$b'$}
\cl{\hspace{-1.5cm}\includegraphics{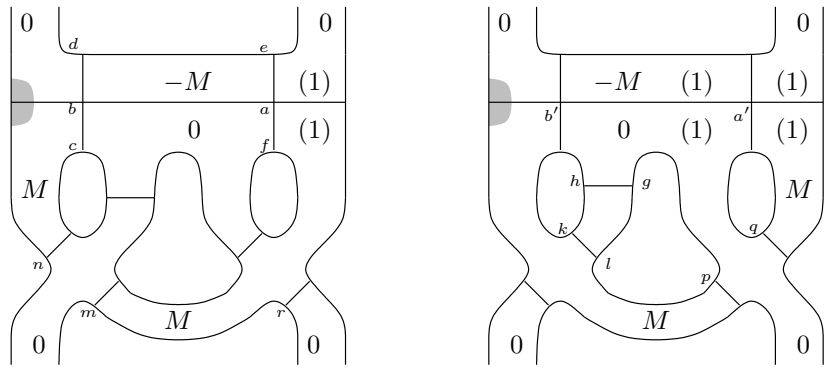}}
\caption{\label{fig:whitehead.4} {The left and the right picture
represent the upper and lower faces  of the surface $C$. The
intersection points of the $\alpha$ and $\beta$ circles are
labelled. There is a set of periodic domains with coefficients
$0,M,-M$ in different regions, as above, for any given number $M$.
There are two distinguished domains connecting the Kauffman
states: (a) The small rectangle with vertices $g,h,k,l$, (b) the
domain whose non-zero coefficients on different domains are
denoted by $(.)$.}}
\end{figure}

If $\{\x,\y\}$ and
$\{\x,\y'\}$ are two meridian Kauffman states, the domain between
 $\y,\y'$ on $\Sigma$ will have  coefficient $0$ on one side
and $M$ on the other side of the meridian $m$. We may complete this
domain to the domain of an actual disk connecting $\{\x,\y\}$
and $\{\x,\y'\}$
by adding the domain on $C$ with the multiplicities shown in \figref{fig:whitehead.4} .

Therefore, any two meridian Kauffman states $\{\x,\y\}$ and
$\{\x,\y'\}$ are in the same $\SpinC$ class.

A similar argument using the periodic domain $\mathcal{D}$ of
 $(\Sigma;\deltas;\{l\}\cup \gammas_0)$ and the above periodic
domain of $C$  in  \figref{fig:whitehead.4},
 shows that any two longitude Kauffman states $\{\x,\y\}$ and
$\{\x,\y'\}$ are also in the same $\SpinC$ class.

To understand the $\SpinC$ grading between meridian and longitude Kauffman states,
the next step is to consider the Kauffman states of the two Heegaard diagrams
\begin{displaymath}
\begin{CD}
H_1=(C;\gamma,\alpha_1,\ldots,\alpha_4;n,\beta_1,\ldots,\beta_4),
H_2=(C;\lambda,\alpha_1,\ldots,\alpha_4;n,\beta_1,\ldots,\beta_4).
\end{CD}
\end{displaymath}
We have named the intersection points of the diagram by letters of the alphabet in
\figref{fig:whitehead.4}.

The
Kauffman states of $H_1$ will be the following list
\begin{displaymath}
\begin{CD}
\x_1=\{n,d,h,p,f\},\\
\x_2=\{n,c,g,r,f\},\\
\x_3=\{n,g,c,q,e\}.
\end{CD}
\end{displaymath}
It is easy to see that $\x_1$ and $\x_2$ are in the same $\SpinC$ class and there
is a disk between them,
which is disjoint from the shaded area, where the handle
is attached to the surface $C$. This disk supports a unique holomorphic representative.
The numbers in parenthesis denote the nonzero coefficients of the disk between these
two Kauffman states.

The Kauffman state $\x_3$ is in a higher $\SpinC$ class. This means that
$\spinc(\x_1)+1=\spinc(\x_2)+1=\spinc(\x_3)$.

The relative Maslov grading of any two Kauffman states (of the
meridian or longitude type) $\{\x,\y\}$ and $\{\x,\y'\}$ is the
same as the relative Maslov index of the states $\y$ and $\y'$. So
the contribution of all the Kauffman states of the form
$\{\x,\y\}$, for a fixed $\x$, is equal (up to a sign) to the
Euler characteristic of $\widehat{HF}(S^3)$ or the Euler
characteristic of $\widehat{HF}(S^3_0(K))$, depending on whether
$\{\x,\y\}$'s are meridian Kauffman states or longitude Kauffman
states, respectively. This implies that the total contribution of
longitude Kauffman states to the Euler characteristic in different
$\SpinC$ structures is zero. For each of $\x_i$ we will get a
contribution equal to $\pm1$. The contribution from $\x_1$ is
cancelled against the contribution from $\x_2$, so the only
$\SpinC$ structure for which the contribution is nonzero, is
$s(\x_3)$. The conclusion is that $s(\{\x_3,\bullet\})=s(\x_3)=0$,
since the Euler characteristic of $\widehat{\HFK}(K_L)$ gives the
symmetrized Alexander polynomial of $K_L$ (which is trivial).
Moreover, as a result of the previous discussion, we have:
$$\spinc(\{\x,\y\})=\spinc(\x),\ \ \ \ \text{for all meridian Kauffman
  states }\{\x,\y\}.$$
Here the right hand side is a $\SpinC$ class associated with the
  Heegaard diagram $H_1$.
So $\spinc(\x_1)+1=\spinc(\x_2)+1=\spinc(\x_3)=0$, and there is no
meridian Kauffman state in the $\SpinC$ class $s=1$.

Now we turn to the Kauffman states of $H_2$.
We will show that some of them naturally cancel against each other. Then we will identify
the remaining ones, and will compute the $\SpinC$-grading of the corresponding longitude Kauffman
states.

There is a small rectangle bounded by the
intersection points $h,g,l,k$ on the lower face of the surface.
For any pair of Kauffman states for the Whitehead double
which are of the form $\{g,k,\bullet\}$ and $\{h,l,\bullet\}$, the domain of the
disk between these two states
is this rectangle, which supports a unique holomorphic representative. There are six of
these pairs. We may cancel them against each other, in the expense that having two Kauffman states
with a disk of the following type connecting them,
we will not be able to argue that there is no
boundary map between the two Kauffman states.
The disks considered above are the ones with negative coefficients in the rectangle.
Since we will not face this situation in
what follows, we simply choose to cancel them against each other. For a more careful
explanation of this method we refer the reader to Rasmussen's paper \cite{Ras}.

Here is a list of the remaining  Kauffman states:
\begin{displaymath}
\begin{CD}
\z_1=\{a,q,g,m,c\},\\
\z_2=\{a',q,g,m,c\},\\
\z_3=\{b,h,m,p,f\},\\
\z_4=\{b',h,m,p,f\}.
\end{CD}
\end{displaymath}
One may check by considering the domains that
\begin{equation}
\spinc(\z_1)=\spinc(\z_2)+1=\spinc(\z_3)+1=\spinc(\z_4)+2.
\end{equation}
In order to see what the absolute grading of these states is, move $\delta_1$ (
the unique $\delta$ curve that intersects $m$) by an isotopy
 to create a pair of intersection points
with $l$. One of them has the property that together with the intersection of $l$ and $m$
and the intersection of $m$ and $\delta_1$, they form the vertices of a small triangle.
Call this point $x_0$ and let $y_0$ be the intersection point of $m$ and $\delta_1$.
If $\y=\{y_0,\bullet\}$ is a Kauffman state for
$(\Sigma;\deltas;\{m\}\cup \gammas_0)$, then $\{x_0,\bullet\}$ will be a
Kauffman state for
$$
(\Sigma;\deltas;\{l\}\cup \gammas_0).$$
There is a domain representing a disk with zero coefficients on $z$ and $w$ which
connects the two Kauffman states $\z_3\cup \{y_0,\bullet\}$ and
$\x_3\cup \{x_0,\bullet\}$. So for any longitude Kauffman state
of the form $\z_i\cup \{\bullet\}$ we may compute the $\SpinC$
grading via the formula:
\begin{equation}
\spinc(\z_1\cup \{\bullet\})-1=\spinc(\z_2\cup \{\bullet\})=
\spinc(\z_3\cup \{\bullet\})=\spinc(\z_4\cup \{\bullet\})+1=0.
\end{equation}
So the only Kauffman states in the $\SpinC$ structure $\spinc=1$ that remain, are those of the form
$\z_1\cup \y$, where $\y$ is some Kauffman state on the
Heegaard diagram $H_2$ (which is a potential Heegaard diagram for the longitude
Floer homology).

Suppose that $\z_1\cup \y$ and $\z_1\cup \y'$
are two Kauffman states in our Heegaard diagram. Since the states
do not differ on $C$, the
only possibility for a domain between the two states is that the coefficients in
all of the regions on $C$ are zero except for the regions where a coefficient equal to
$\pm M$ is assigned as in \figref{fig:whitehead.4}.
In any such domain, there are regions with both
$M$ and $-M$ as coefficients. Furthermore, these domains do not use the small rectangle
considered before. So the only case where there is potentially a boundary map from
$\z_1\cup \y$ to $\z_1\cup \y'$ is when we have $M=0$.

In this case the four regions around the connecting handle will get coefficients equal to
zero. This means that the disk is completely supported on $\Sigma$. Furthermore
if we put two marked points $z'$ and $w'$ on the two sides of $l$ at the intersection of $l$ with
$\beta_1$, the above discussion shows that the domains of the disks between these points will have
zero coefficients in the regions associated with $z'$ and $w'$.

So, the  disks that contribute to the boundary operator are in 1-1 correspondence
with the disks between $\y$ and $\y'$ in the hat theory assigned to
the Heegaard diagram
$$
(\Sigma,\deltas;\{l\}\cup \gammas_0;z',w').$$
The above discussion shows that the generators and all the boundary maps in the
Ozsv\'ath-Szab\'o Floer homology of $K_L$ in $\SpinC$ structure $s=1$ are exactly
the same as those appearing in
$$
\widehat{\CFL}(K)=\bigoplus_{\ui\in \mathbb{Z}+\frac{1}{2}}\widehat{\CFL}(K,\ui).
$$
The $\SpinC$ grading of $\widehat{\CFL}(K)$, and that of the homology groups
$\widehat{\HFL}(K)$ are forgotten when we compute the Ozsv\'ath-Szab\'o
Floer homology of the Whitehead double, and the isomorphism is an isomorphism
of groups (relatively) graded by the Maslov index.

We have proved the following theorem:

\begin{thm}{Let $K_L$ denote the Whitehead double of a knot $K$ in $S^3$. The
Ozsv\'ath-Szab\'o Floer homology groups $\widehat{\HFK}(K_L,\pm 1)$ are isomorphic to the
group $\widehat{\HFL}(K)=\bigoplus_{\ui\in
  \mathbb{Z}+\frac{1}{2}}\widehat{\HFL}(K,
\ui)$ as (relatively)
$\mathbb{Z}$-graded abelian groups with the (relative) grading on both sides coming from the
Maslov grading.}
\end{thm}
As a corollary of this theorem and the results of the previous section we have:

\begin{cor}{ Let $K=T(2,2n+1)$ denote the $(2,2n+1)$ torus knot and
let $K_L$ be the Whitehead double of $K$. Then the Ozsv\'ath-Szab\'o Floer homology groups
$\widehat{\HFK}(K_L,+1)$ in different (relative) Maslov gradings are described as follows:
\begin{displaymath}
\begin{array}{cccccc}
\mu: & n& n-2 & \ldots & -n+2 &-n+1\\
\widehat{\HFK}:&
\Z\oplus \Z& \Z \oplus \Z &\ldots & \Z\oplus \Z & \bigoplus_{i=1}^{2n}\Z\\
\end{array}
\end{displaymath}
}\end{cor}

\begin{remark}
The longitude Floer homology may be defined for a knot in a
three-manifold $Y$, and as such, it enjoys very nice surgery formulas.
We postpone a discussion of these subjects to a future paper.
\end{remark}

\Addresses\recd

\begin{thebibliography}{99}
\bibitem{Burde}
\textbf{G Burde}, \textbf{H Zieschang}, \emph{Knots}, de Gruyter Studies in
  Mathematics 5, Walter de Gruyter \& Co. Berlin (2003) \MR{1959408}

\bibitem{OS-knot}  \textbf{P Ozsv{\'a}th}, \textbf{Z Szab{\'o}}, 
\emph{Holomorphic disks and knot invariants}, to appear in Advances in Math,
\arxiv{math.GT/0209056}

\bibitem{OS-alternating}
\textbf{P Ozsv{\'a}th}, \textbf{Z Szab{\'o}}, \emph{Heegaard {F}loer homology
  and alternating knots}, \gtref7{2003}{6}{225}{254} \MR{1988285}

\bibitem{OS-3m1} \textbf{P Ozsv{\'a}th}, \textbf{Z Szab{\'o}},
\emph{Holomorphic disks and topological invariants for closed
three-manifolds}, to appear in Annals of Math. \arxiv{math.SG/0101206}

\bibitem{OS-3m2} \textbf{P Ozsv{\'a}th}, \textbf{Z Szab{\'o}},
\emph{Holomorphic disks and three-manifold invariants: properties and
applications}, to appear in Annals of Math. \arxiv{math.SG/0105202}

\bibitem{OS-fibered} \textbf{P Ozsv{\'a}th}, \textbf{Z Szab{\'o}}, 
\emph{Heegaard Floer homologies and contact structures},
\arxiv{math.SG/0210127}

\bibitem{OS-4genus}
\textbf{P Ozsv{\'a}th}, \textbf{Z Szab{\'o}}, \emph{Knot {F}loer homology and
  the four-ball genus}, \gtref7{2003}{17}{615}{639} \MR{2026543}

\bibitem{OS-genus}
\textbf{P Ozsv{\'a}th}, \textbf{Z Szab{\'o}}, \emph{Holomorphic disks and genus
  bounds}, \gtref8{2004}{8}{311}{334} \MR{2023281}

\bibitem{Ras}
\textbf{J\,A Rasmussen}, \emph{Floer homology of surgeries on two-bridge
  knots}, \agtref2{2002}{32}{757}{789} \MR{1928176}

\bibitem{Ras2} textbf{J\,A Rasmussen}, \emph{Floer homology and knot
complements}, PhD thesis, Harvard Univ. \arxiv{math.GT/0306378}

\bibitem{Rud}
\textbf{L Rudolph}, \emph{The slice genus and the {T}hurston-{B}ennequin
  invariant of a knot}, Proc. Amer. Math. Soc. 125 (1997) 3049--3050
  \MR{1443854}

\end{thebibliography}
\end{document}